%This is an AMS-TeX file
\input amstex
%\loadbold

\documentstyle{amsppt}
\magnification=\magstep1
%\magnification\magstep2
\document     
\baselineskip=12pt
\parindent=10pt
\define\e{\varepsilon}

\define\ztil{\tilde{z}}

\define\etatil{\tilde{\eta}}
\define\ubar{\bar{u}}
\define\Xbar{\bar{X}}
\define\zbar{\bar{z}}
\define\xbar{\bar{x}}
\define\leb{\operatorname{Leb}}
\define\lo{\operatorname{LO}}
\define\essinf{\operatornamewithlimits{ess\, inf}}
\define\esssup{\operatornamewithlimits{ess\, sup}}
\define\Se{Sepp\"al\"ainen (1999)}

\define\aihp{Ann.\ Inst.\ H.\ Poincar\'e Probab.\ Statist.}
\define\aihpnl{Ann.\ Inst.\ H.\ Poincar\'e Anal.\ Non Lin\'eaire}
\define\jsp{J.\ Statist.\ Phys.}

\define\ap{Ann.\ Probab.}
\define\mprf{Markov Process.\ Related Fields}
\define\ptrf{Probab.\ Theory Related Fields}

\define\zwvg{Z.\ Wahrsch.\ Verw.\ Gebiete}

\centerline{\bf Second class particles
as microscopic characteristics in }
\centerline{\bf  totally asymmetric nearest-neighbor $K$-exclusion
processes}

\hbox{}

\centerline{2000}

\hbox{}

\centerline{Timo Sepp\"al\"ainen 
\footnote""{ Research partially supported by NSF grant DMS-9801085. }}
\hbox{}
\centerline{Department of Mathematics}
\centerline{Iowa State University}
\centerline{Ames, Iowa 50011, USA}
\centerline{seppalai\@iastate.edu}

\hbox{}

\hbox{}

\flushpar
{\it Summary.}
We study aspects of the hydrodynamics 
of one-dimensional totally asymmetric $K$-exclusion, building on 
the hydrodynamic limit of \Se.  We prove that 
the weak solution chosen by the particle system is 
the unique one with maximal current past any fixed 
location. A uniqueness result is needed
because we can
prove neither differentiability nor strict concavity 
of the flux function, so we cannot
use the Lax-Oleinik formula or jump conditions to 
define entropy solutions. Next we prove laws
of large numbers for a second class particle 
in $K$-exclusion. The macroscopic trajectories of
second class particles are characteristics and shocks of the
conservation law for the particle density. 
In particular, we extend to $K$-exclusion 
Ferrari's result
that the second class particle 
follows a macroscopic shock in the Riemann solution.
The technical novelty of the proofs is  a variational
representation for the position of a second class
particle, in the context of the variational coupling
method.

\hfil

\flushpar
Mathematics Subject Classification: Primary
 60K35, Secondary  82C22 

%\hfil

\flushpar
Key words: Exclusion process, second class particle, 
characteristic, 
 hydrodynamic limit, variational coupling method 

%\hfil

\flushpar
Short Title: Second class particles in $K$-exclusion

\vfill
\break

\head 1. Introduction \endhead

This paper continues the study of the hydrodynamics
of totally asymmetric nearest-neighbor $K$-exclusion
processes. These are interacting systems 
of 
indistinguishable particles  on the one-dimensional 
integer lattice $\bold Z$, subject to the restriction 
that each site contain at most $K$ particles. 
Independently at each site and at constant exponential 
rate one, one particle jumps  to the next
site on the right, provided the $K$-exclusion rule is not
violated. 

A hydrodynamic limit
for this process was proved in \Se. 
This is a law of large numbers 
according to which  the empirical particle density
converges 
to a weak solution of a nonlinear conservation law, 
under the appropriate scaling of space
and time. 

The present paper begins with the question of
uniqueness of the weak solution chosen by the 
particle system.   The criterion we give 
is maximization of current across any fixed position.
 This issue arises because 
we do not know whether the flux function of the 
conservation law  is strictly concave or differentiable.
 Consequently we cannot use
the Lax-Oleinik formula to define an entropy solution,
or one-sided bounds on jumps as an entropy criterion. 

We address 
this difficulty because it is closely
 tied to the main topic of the paper,
the dynamics of a second class
particle in the hydrodynamic scaling. 
As we shall see, our incomplete understanding 
of the macroscopic equation constrains the 
macroscopic description of the second class
particle. 

A {\it second class
particle} is an extra particle in the system that 
 yields right of way to the regular, first class
particles. The second class particle can jump only
if there is no first class particle at its site who
can jump, and a second class particle has to jump 
backward if first class particles fill its site. 

Second class particles in asymmetric exclusion
processes have been studied 
 for several reasons. In a certain
technically precise sense, the second class particle 
marks  the location of a shock at the 
microscopic particle level. These ideas were 
developed by Ferrari, Kipnis and Saada (1991),
Ferrari (1992), and 
Ferrari and Fontes (1994a).  At this same 
time, Ferrari (1992) discovered that in the 
hydrodynamic limit a second class particle follows
the characteristics of the macroscopic equation. 
The connection between second class particles and 
characteristics was further elucidated by Rezakhanlou
(1995) and Ferrari and Kipnis (1995). 
The second class particle also enters naturally
in the study of the fluctuations of the 
current, as in Ferrari and Fontes (1994b). 

Another strain of work, technically related to the 
investigation of microscopic shocks, concerns 
the invariant distributions of the exclusion process
as seen from a second class particle, and
the invariant distributions of the 
two-species process that has positive densities
of both first and second class particles. This
program culminated in the complete description 
of these measures by 
Ferrari, Fontes and Kohayakawa (1994).
These and other results, such as the central limit
behavior of the second class particle, are discussed 
in Part III of Liggett (1999) and in the notes
of Chapter 9 in Kipnis and Landim (1999). 

The above results are for systems whose invariant
distributions are explicitly known, namely the 
exclusion and zero-range processes. And 
the proofs typically use precise and deep knowledge
about those invariant measures. 

In the present paper we prove laws of large numbers
for a single second class particle in the 
hydrodynamic scaling of totally asymmetric $K$-exclusion. 
Essentially no information 
is currently available about the invariant distributions
of these asymmetric $K$-exclusion processes. 

Of the papers mentioned above, 
 Ferrari (1992) and Rezakhanlou (1995)
are the closest predecessors of our paper. 
Both derived laws of large numbers 
for the second class
particle. The limiting trajectory is a characteristic or a 
 shock of the conservation law of the 
particle density. 
 Ferrari proved 
strong laws  in an asymmetric simple exclusion 
process, mainly in these settings:
 the initial particle 
distribution is either an i.i.d.\ product measure,
or a product measure
that represents a shock profile that is constant
on either side of the origin (this is known as
the Riemann problem).  
Ferrari's proofs rely on 
couplings of processes of several species of
particles, and knowledge of invariant
distributions. 

Rezakhanlou (1995) proved $L^1$-laws of large 
numbers for the second class
particle in asymmetric simple exclusion and 
in a totally asymmetric zero-range process. He
assumed that 
 initial distributions are  product measures, and 
the initial  
macroscopic profiles are bounded integrable 
functions on $\bold R$, 
with some additional technical assumptions. 
His proof derives the 
limit of a single second class
particle from the limit of a macroscopically
visible interval  of second class
particles, and the latter he gets
from the basic hydrodynamic limits
for exclusion and zero-range processes, from 
Rezakhanlou (1991). The invariant distributions of the 
processes enter his proofs through the macroscopic
equations. These   are known explicitly
precisely because the invariant distributions are known. 
Rezakhanlou's approach would prove some 
result for $K$-exclusion, since one can 
use the hydrodynamic limit of \Se\ instead of 
Rezakhanlou (1991).  

Our approach is  different from the earlier ones.
We use  a representation of the 
 second class particle in the variational
coupling method. To preview 
 this, let us switch from the 
exclusion process to the associated marching soldiers 
process (in other words, consider the interface
process instead of the increments process). 
In this approach a process $z(i,t)$ from a general
initial interface is coupled 
with auxiliary processes $w^k(i,t)$ that evolve from 
translated  wedge-shaped initial profiles. The 
coupling is arranged so that $z(i,t)$ is
 the envelope of
the auxiliary processes:
$$z(i,t)=\sup_{k\in\bold Z} w^k(i-k,t).
\tag 1.1
$$
Let $X(t)$ denote the position of a second class particle 
added to the  exclusion process  
$\eta(i,t)=z(i,t)-z(i-1,t)$.  
The new tool is this variational expression: 
$$X(t)=\inf\{i: \text{$z(i,t)=w^k(i-k,t)$ for 
some $k\ge X(0)$} \}.
\tag 1.2
$$
Variational formula (1.1) is a
microscopic analogue of the Hopf-Lax formula 
of viscosity solutions of Hamilton-Jacobi
equations. In this same vein, (1.2) corresponds
to the definition of forward characteristics
 through the Hopf-Lax formula. This  
justifies regarding  the second class particle as
a microscopic analogue of the characteristic,
as expressed in the title of the paper.

Variational formulas (1.1)--(1.2) permit us to work
directly on the paths of the process,
and we do not need a particular form
for the  initial particle distributions. 
We need to assume that the initial particle profile 
and the initial position of the second class
particle both converge under the 
appropriate space scaling. The result is
a law of large numbers for the second class
particle.  The theorems give either convergence
in probability or almost surely, depending on which
kind is assumed initially. 
The limit is a 
characteristic or a shock
of the macrosopic partial differential
equation.

Here we run into the difficulty mentioned above. For
$2\le K<\infty$ (cases other than exclusion or 
zero-range) we do not have sufficient control of the
regularity of the flux $f$ of the conservation law. 
For this reason we cannot assert that macroscopically 
the second class particle follows the solution of
$\dot x=f'(\rho(x,t))$, where $\rho(x,t)$ is the 
macroscopic particle density. Instead, we use the
Hopf-Lax formula to define the macroscopic trajectory
of the second class particle.

If the macroscopic
characteristics are unique (there are no
rarefaction fans) we do not need the assumption 
of an initial law of large numbers for the 
second class particle. In this case the second class
particle follows, in the hydrodynamic scale,
 the macroscopically defined characteristic
that emanates from its random initial position. 

We improve earlier results  
in two ways: by broadening the class of processes
beyond exclusion and zero-range, and by permitting
more general initial distributions. 
In one important sense the results of Ferrari and 
Rezakhanlou remain more general than ours: their results
cover an asymmetric simple exclusion process where particles
can jump both right and left, while our results are 
restricted to {\it totally} asymmetric processes
where particles jump only to the right. This is 
because a way to apply the variational
coupling approach without total asymmetry has yet to be  
discovered. 

\hbox{}

{\it Overview of the paper.}
In Section 2 we discuss the results. 
Section 2.1  reviews the hydrodynamics of totally
asymmetric $K$-exclusion
and zero-range processes (Theorem 1).  
We state a uniqueness 
criterion for the weak solution chosen by the particle 
system, valid for a concave flux that is not necessarily
differentiable or strictly concave (Theorem 2).    
Section 2.2 states the laws of large 
numbers for the second class particle (Theorems
3 and 4). Section 2.3 
comments on the various assumptions and conclusions. 

In Section 2.4 we
 discuss a generalization of the 
$K$-exclusion and zero-range processes to which 
our results apply. In this model
 particles may jump in 
batches of varying sizes
from one site  to the next. Such a generalization 
is natural from a queueing perspective, where the
particles are customers moving through
a sequence of service stations. 

Some standard tools of conservation
laws are not at our disposal for lack of
 information about the regularity of  the 
flux $f$ of the conservation law. 
  This problem is clarified in
Section 3.  Without strict concavity of $f$  
the Lax-Oleinik formula
cannot be used, and the correct
 weak solution can have 
 shocks prohibited by entropy criteria (Section 3.1). 
Without
differentiability of $f$ we loose some good properties
of characteristics (Section 3.2). In section 3.3 we prove Theorem 2,
the uniqueness result.

Section 4 reviews the essentials of the variational
coupling approach to exclusion and zero-range processes from
Sepp\"al\"ainen (1998a,b,1999), and develops the 
formula for the second class particle. Section 5
proves the laws of large numbers for 
the second class particle. The proofs are an 
interplay of the particle level properties 
from Section 4 and the properties of the 
macroscopic characteristics. 

{\it Notation.} The set of natural numbers is 
$\bold N=\{1,2,3,\ldots\}$, and the set of nonnegative
integers $\bold Z_+=\{0,1,2,3,\ldots\}$.

\head 2. Results \endhead

\subhead 2.1. Hydrodynamic limits and entropy solutions
\endsubhead
We start by reviewing the hydrodynamic limit of
totally asymmetric $K$-exclusion. 
The parameter $K$ is a fixed
 positive finite or infinite integer $K\in\bold N\cup\{\infty\}$.
Each process consists of indistinguishable particles 
that move on the 
one-dimensional integer lattice $\bold Z$, 
subject to the restriction that a site 
$i\in\bold Z$ can be occupied
by at most $K$ particles.  If $K=\infty$ this 
constraint is absent. The state space of the process is 
$\Cal S=\{0,1,2,\ldots, K\}^{\bold Z}$ for finite $K$, 
and $\Cal S={\bold Z}_+^{\bold Z}$ if $K=\infty$. 
 A generic
element of $\Cal S$ is
 $\eta =(\eta(i):i\in\bold Z)$, where  $\eta(i)$ 
denotes the number
of particles present at site $i$.

The dynamics is such that at exponential rate 1, 
one particle moves from site $i$ to $i+1$,
provided there is  at least one  particle at $i$ 
and at most $K-1$ particles at $i+1$. 
 The jump attempts happen independently 
at all sites $i$. 

 To realize this dynamics, let $\{D_i:i\in\bold Z\}$
denote a collection of mutually independent rate 1  
Poisson point processes on the time axis
$(0,\infty)$. 
We usually think of $D_i$ as the random set of 
jump times (or {\it epochs}, to borrow Feller's term) rather
than as the corresponding random step function.  
At each epoch of $D_i$ a  jump
from 
site $i$ to site $i+1$ is attempted, and 
successfully executed if the rules permit. 
The state of the process at time $t$ 
is the configuration
$\eta(t) =(\eta(i,t):i\in\bold Z)$ of occupation numbers.

The dynamics can be
 represented by the generator $L$  that
acts on bounded
cylinder functions $f$ on $\Cal S$:
$$Lf(\eta)=\sum_{i\in\bold Z} \bold 1\{\eta(i)\ge 1,\eta(i+1)\le K-1\}
[f(\eta^{i,i+1})-f(\eta)].
\tag 2.1
$$
Here $\eta^{i,i+1}$ is the configuration that results from the
 jump of a single particle  from site $i$ to site $i+1$: 
$$ 
\eta^{i,i+1}(j)=
\cases \eta(i)-1, &j=i\\
\eta(i+1)+1, &j=i+1\\
\eta(j), &j\ne i,i+1.
\endcases
$$

{\it Special cases.} 
(a) When $K=1$ this process is known as the totally asymmetric simple
exclusion process (TASEP).  The Bernoulli probability measures
 $\{\nu^\rho\}$ on $\Cal S$, 
indexed by density $\rho\in[0,1]$, are invariant for
the process. These are defined by the condition
$$\nu^\rho\{ \eta(i_1)=\eta(i_2)=\cdots=\eta(i_m)=1\}=\rho^m
$$
for any finite set $\{i_1,\ldots,i_m\}$ of sites. 

(b) If  $K=\infty$ this process is a member of the family
of totally asymmetric zero-range processes (TAZRP).  
The i.i.d.\ geometric probability measures
 $\{\nu^\rho\}$, $\rho\in[0,\infty)$,
 are invariant. These are given 
by 
$$\nu^\rho\{ \eta(i_1)=k_1, \eta(i_2)=k_2,
\ldots,\eta(i_m)=k_m\}=\prod_{s=1}^m \frac1{1+\rho}
\left(\frac\rho{1+\rho}\right)^{k_s}
\tag 2.2
$$
for distinct sites $i_1,\ldots,i_m$ and arbitrary
integers $k_1,\ldots,k_m\in\bold Z_+$.

(c) For  $2\le K<\infty$ we 
 call this the totally asymmetric $K$-exclusion
process. For this case the existence of spatially ergodic
invariant measures for all densities has not been proved.

\hbox{}

The basic result about the large-scale behavior of these
processes is the following hydrodynamic scaling limit. 
Let 
 $\rho_0$ be  a bounded measurable function on $\bold R$ such 
that 
$0\le \rho_0(x)\le K.
$
Assume given 
 a probability space on which are defined a sequence of
 initial configurations $(\eta_n(i,0):i\in\bold Z)$, 
 $n\in\bold N$, and the Poisson jump time processes
$\{D_i\}$.  The processes $\{\eta_n(\cdot)\}$ are
constructed on this probability space by the familiar
graphical method.
 Assume that 
at time zero the following strong law of large numbers 
is valid: 
$$
\text{For all $-\infty<a<b<\infty$, }\;
\lim_{n\to\infty}\frac1n\sum_{i=[na]+1}^{[nb]}\eta_n(i,0)
=\int_a^b\rho_0(x)dx\ \ \text{a.s.}
\tag 2.3
$$ 

\proclaim{Theorem 1}  Under
assumption {\rm (2.3)} a strong law of large numbers 
holds for all $t>0$ and all $-\infty<a<b<\infty$:
$$
\lim_{n\to\infty}\frac1n\sum_{i=[na]+1}^{[nb]}\eta_n(i,nt)
=\int_a^b\rho(x,t)dx\ \ \text{a.s.}
\tag 2.4
$$ 
The deterministic limiting
 density $\rho(x,t)$ is a weak 
solution of the Cauchy problem 
$$\rho_t+f(\rho)_x=0\,,\quad \rho(x,0)=\rho_0(x).
\tag 2.5
$$
The flux function  $f$ of the conservation law
{\rm (2.5)} is a fixed, deterministic 
concave function on $[0,K]$ that depends on $K$.  
\endproclaim

The limit (2.3) is the only assumption needed in the 
following strong sense:  The joint distribution of the 
initially defined random variables 
$  \{(\eta_n(i,0):i\in\bold Z): n\in\bold N\}$, 
and $\{D_i:i\in\bold Z\}$
may be arbitrary. As long as the marginal distribution 
of $\eta_n(i,0)$ satisfies (2.3), and marginally 
$\{D_i\}$ are i.i.d.\ rate 1 Poisson point processes,
Theorem 1 is valid. 

A second point is that if (2.3) holds in probability,
then the conclusion (2.5) also holds in probability.

 For exclusion ($K=1$) and zero-range 
($K=\infty$) Theorem 1 has gone through many stages
of generalization, beginning with the seminal work
of Rost (1981). The entropy technique of  
Rezakhanlou (1991) proves it for multidimensional 
exclusion and zero-range processes, without requiring 
total asymmetry. The variational coupling method 
of Sepp\"al\"ainen (1998a,b)
permits the most general initial distributions 
for the one-dimensional, totally asymmetric case. 
For $2\le K<\infty$ Theorem 1 was first proved in
\Se.

Let us write $f_K$ for the flux function  $f$  in (2.5) 
when we 
emphasize dependence  on $K$. 
$f_K$  is known in cases
where invariant distributions are explicitly known:  
for exclusion
$f_1(\rho)=\rho(1-\rho)$ for $0\le \rho\le 1$, and 
for zero-range $f_\infty(\rho)=
{\rho}/{(1+\rho)}
$ 
for $0\le \rho <\infty$. 
For $1<K<\infty$ we know the concavity
and continuity of $f_K$, the symmetry 
$f_K(\rho)=f_K(K-\rho)$, the monotonicity 
$f_K(\rho)\le f_{K+1}(\rho)$, and these bounds:
$$
\text{$\min\{\rho(1-\rho)\,,\,1/4\}\le f_K(\rho)
\le \frac{\rho}{1+\rho}$ for $0\le\rho\le K/2$.}
\tag 2.6
$$

Currently we have no information about the
regularity of $f$ for $1<K<\infty$. As suggested by $f_1$ and $f_\infty$,
we would expect $f$ to be $C^1$ and 
strictly concave, in other words, for its graph
to have no
corners or line segments. This lack of information
about $f$ 
creates  a  problem because the common definitions of 
entropy solutions of (2.5) are valid only for strictly
convex or concave flux
functions. Thus we cannot apply 
 standard entropy criteria
such as one-sided bounds on jumps or the Lax-Oleinik
formula. We show this in Section 3. 

The next task is to complement Theorem 1 with a 
uniqueness criterion. 
First we explain how the solution $\rho(x,t)$ 
that the particle system chooses in Theorem 1
is defined.
 In the next discussion 
assume first $K<\infty$. Afterwards we comment on the 
$K=\infty$ case.  

 Define an 
antiderivative $u_0$ of $\rho_0$ by 
$$\text{$u_0(0)=0$  and }\quad
u_0(b)-u_0(a)=\int_a^b \rho_0(x)dx\quad\text{for all $a<b$.}
\tag 2.7
$$
Let $g$ be the nonincreasing, nonnegative convex function on $\bold R$
defined by 
$$g(x)=\sup_{0\le\rho\le K}\{ f(\rho)-x\rho\}.
\tag 2.8
$$
Symmetry of $f$ and bounds (2.6) imply that 
$f'(0)=1$ and $f'(K)=-1$, the existence of the 
derivatives being part of the conclusion. It follows
that $g$  
satisfies 
$$\text{$g(x)=0$ for  $x\ge 1$ and $g(x)=-Kx$ for $x\le -1$.}
\tag 2.9
$$
Then for $x\in\bold R$ set $u(x,0)=u_0(x)$, and for $t>0$ 
$$u(x,t)=\sup_{y\in\bold R}\biggl\{ u_0(y)-tg\biggl(\frac{x-y}t\biggr)
\biggr\}.
\tag 2.10
$$
The supremum is attained at some $y\in[x-t,x+t]$. The 
function $u$ is uniformly Lipschitz  on $\bold R\times[0,\infty)$,
nonincreasing in  $t$ and  nondecreasing in $x$. 
(2.10) is known as the {\it Hopf-Lax formula}.
It defines  $u(x,t)$
as a {\it viscosity solution} of the  Hamilton-Jacobi
equation 
$$u_t+f(u_x)=0\,,\qquad u(x,0)=u_0(x) 
\tag 2.11
$$
[Ch.\ 10 in Evans (1998)]. The viscosity solution 
is unique by Theorem 2.1 of Ishii (1984) that applies
to (2.11) because $f$ is continuous and 
 $u$ is uniformly continuous. 

If $K=\infty$ we know that $g(x)=(1-\sqrt{x}\,)^2$ for 
$0\le x\le 1$, $g\equiv 0$ on $(1,\infty)$ and we can 
take $g\equiv \infty$ on $(-\infty, 0)$. The supremum in 
(2.10) is attained at some $y\in[x-t,x]$. 
Since $\rho_0$ is assumed bounded, 
$u$ is again Lipschitz  on $\bold R\times[0,\infty)$.

Let 
$$\rho(x,t)=\frac{\partial}{\partial x}u(x,t)
\tag 2.12
$$ 
be the a.e.\ defined $x$-derivative. Recall this 
standard definition: 
 a locally bounded measurable 
function $\lambda(x,t)$ is a {\it weak solution}
of (2.5) if for all 
compactly supported and continuously differentiable
test functions $\phi$ on $\bold R\times[0,\infty)$, 
$$\int_0^\infty\int_{\bold R} [\lambda(x,t) \phi_t(x,t)+
f(\lambda(x,t))\phi_x(x,t)]dxdt
+\int_{\bold R} \rho_0(x)\phi(x,0) dx=0.
\tag 2.13
$$
By the proof
of Thm.\ 3.4.2 in Evans (1998), $\rho(x,t)$ is a 
weak solution of (2.5) in the integral sense (2.13).
The solution $\rho(x,t)$ defined by (2.12) is the one chosen 
by the particle system in Theorem 1. 

In cases
$K=1$ and $K=\infty$ we know that the flux $f$ 
is differentiable and strictly concave. Then 
$\rho(x,t)$ can be equivalently defined by the 
Lax-Oleinik formula, 
and we can call $\rho(x,t)$
the ``unique entropy solution.'' For 
$2\le K<\infty$  we can only define $\rho(x,t)$
through (2.12). In Sections 3.1--3.2
we show that without further knowledge about $f$
this solution cannot be defined by the Lax-Oleinik
formula, and we cannot rule out  discontinuities
of the wrong kind that violate entropy conditions. 

To complement the characterization of $\rho(x,t)$ through
  the viscosity solution
of (2.11), we give  a uniqueness criterion 
for which strict concavity
or differentiability of $f$ are not needed.

\proclaim{Theorem 2} Among the weak solutions of {\rm (2.5)}, 
$\rho(x,t)$ of {\rm (2.12)} is characterized by
maximal current over time. More precisely, 
suppose $\lambda(x,t)$ is a nonnegative, locally bounded
measurable function on $\bold R\times(0,\infty)$
that satisfies the integral
criterion {\rm (2.13)}. Fix  $t>0$. Then for
almost all $x\in\bold R$, 
$$\int_0^t f(\lambda(x,s))ds\le \int_0^t f(\rho(x,s))ds.
\tag 2.14
$$
Conversely, if equality holds in {\rm (2.14)} for
almost all $(x,t)$, then $\lambda(x,t)=\rho(x,t)$ 
for almost all $(x,t)$. 
\endproclaim

Thus the particle system chooses the weak solution that
transports the most material. Theorem 2 is proved 
in Section 3.3. 
The system maximizes
the current because $f$ is concave. With a convex
flux the current would be minimized, as is the case with
the stick model version of Hammersley's process, 
see Sepp\"al\"ainen (1996).

\subhead 2.2. Second class particles and characteristics
\endsubhead 
In this section we state the laws of large
numbers for the second class particle.
The 
representation of the second class particle in the 
variational coupling is in 
Section 4.2. 

Let $X(t)$
denote  the position of a single  second class
particle added to the process $\eta(t)$. 
 $X(t)$ is rigorously defined as follows.
Assume given a probability space on which are defined
the following random variables: 
the Poisson jump time processes $\{D_i\}$, an 
initial particle configuration $(\eta(i,0):i\in\bold Z)$, and 
an initial location $X(0)$. 
Define a second initial configuration $(\etatil(i,0):i\in\bold Z)$
that differs from $\eta$ only at site $X(0)$: 
$\etatil(X(0),0)=\eta(X(0),0)+1$ and $\etatil(i,0)=\eta(i,0)$
for $i\ne X(0)$. [Of course we must have 
$\eta(X(0),0)\le K-1$ a.s.\ for this to be possible.]

Run the processes  $\eta(t)$ and $\etatil(t)$ 
so that they read the same Poisson jump time processes $\{D_i\}$.
This is known as the basic coupling of 
$\eta(t)$ and $\etatil(t)$. 
There is always a unique site $X(t)$ at which the two 
processes differ: 
$\etatil(X(t),t)=\eta(X(t),t)+1$ and $\etatil(i,t)=\eta(i,t)$
for $i\ne X(t)$. This defines the position  
$X(t)$  of the second class
particle.

The way in which the second class particle yields 
to first class particles becomes evident when we
consider the two types of transitions that can happen
 to it. Let $X(\tau-)=i$.  

(i) Suppose $\tau$ is a jump epoch for 
Poisson process $D_i$. The second class particle
can jump from $i$ to $i+1$ if there is space at site $i+1$,
but 
only if there is no  first class particle to jump.
In other words, if  $\eta(i,\tau-)=0$
and $\eta(i+1,\tau-)< K$, then 
$X(\tau)=i+1$. 
Otherwise the second class particle does not move
 and $X(\tau)=X(\tau-)=i$. 

(ii) Suppose $\tau$ is a jump epoch for 
$D_{i-1}$. The second class particle
jumps from $i$ to $i-1$
if first class particles fill site $i$.
In other words, if $\eta(i-1,\tau-)\ge 1$ and  $\eta(i,\tau-)=K-1$,
then  $X(\tau)=i-1$.

If $K=\infty$ the second class particle never jumps
left and only the first kind of transition happens to it. 

\hbox{}

 The definition 
in terms of $\eta$ and $\etatil$ generalizes 
naturally to a two-species process of first and second class
particles.   But in this paper we consider only a single second class
particle added to the process.

Next we discuss the 
characteristics of the conservation law (2.5). 
The lack of information about the regularity of $f$
becomes a problem again.
The well-known description of characteristics as
generalized solutions of the ordinary differential 
equation $dx/dt=f'(\rho(x,t))$ 
[see Dafermos (1977)] does not work now. 
As Rezakhanlou (1995), we use the Hopf-Lax formula 
to  define characteristics.  

Once $u(x,t)$ is defined by (2.10), 
 a semigroup property is in force: for all $0\le s<t$
and $x\in\bold R$, 
$$u(x,t)=\sup_{y\in \bold R}
\biggl\{ u(y,s)-(t-s)g\biggl(\frac{x-y}{t-s}\biggr)
\biggr\}.
\tag 2.15
$$
The supremum is assumed at some point 
$y\in [x-(t-s),x+(t-s)]$. (If $K=\infty$
this range would actually be $[x-(t-s),x]$.) By continuity 
the supremum in (2.15) is attained 
at the least and 
largest maximizers, denoted by
$$y^-(x;s,t)=\inf\left\{ y\ge x-(t-s): u(x,t)=u(y,s)-(t-s)
g\biggl(\frac{x-y}{t-s}\biggr) \right\}
\tag 2.16
$$
and 
$$y^+(x;s,t)=\sup\left\{ y\le x+(t-s): u(x,t)=u(y,s)-(t-s)
g\biggl(\frac{x-y}{t-s}\biggr) \right\}. 
\tag 2.17
$$
Abbreviate 
$y^\pm(x,t)=y^\pm(x;0,t)$. As functions of $x$ for fixed $0\le s<t$, 
$y^+(x;s,t)$ is right-continuous, 
$y^-(x;s,t)$ is left-continuous,  both are nondecreasing, 
and trivially from the definition  
$
\lim_{x\to-\infty}y^\pm(x;s,t)= -\infty
$ and $
\lim_{x\to\infty}y^\pm(x;s,t)= \infty.
$
See comment 2.3.1 below for justification of the 
restriction $x-(t-s)\le y\le x+(t-s)$ in 
the definitions of $y^\pm(x;s,t)$.

The minimal and maximal forward characteristics are
defined for $b\in\bold R$, $0\le s<t$, as
$$x^-(b;s,t)=\inf\{ x: y^+(x;s,t)\ge b\}
\tag 2.18
$$
and 
$$x^+(b;s,t)=\sup\{  x: y^-(x;s,t)\le b\}. 
\tag 2.19
$$
Immediate properties are these: 
$$b-(t-s)\le x^-(b;s,t)\le x^+(b;s,t)\le b+(t-s).
$$
As functions of $b$, 
 $x^+(b;s,t)$ is right-continuous, 
$x^-(b;s,t)$ is left-continuous, and both are nondecreasing. 
For the $s=0$ case abbreviate $x^\pm(b,t)=x^\pm(b;0,t)$.
As functions of $t$, $x^\pm(b;s,t)$ are Lipschitz and 
satisfy the initial condition $x^\pm(b;s,s)=b$. 
When characteristics are unique we write
$x(b;s,t)=x^\pm(b;s,t)$.

Now consider the setting of the hydrodynamic limit
Theorem 1, with a sequence of processes $\eta_n$
that satisfy the initial limit (2.3). 
Let $X_n(t)$ denote the position of a second class
particle in the process $\eta_n$. The 
initial position $X_n(0)$ may be deterministic or
random, and 
may depend in an arbitrary fashion
on the initial particle configurations 
$(\eta_n(i,0))$ and the Poisson processes
$\{D_i\}$. For the first result assume a law of large 
numbers at time zero: for a deterministic point $b\in\bold R$
$$
\lim_{n\to\infty} \frac{X_n(0)}n = b \quad\text{a.s.} 
\tag 2.20
$$

\proclaim{Theorem 3} 
Under assumptions {\rm (2.3)} and {\rm (2.20)}, 
for any fixed $t>0$, 
$$
x^-(b,t)\le \liminf_{n\to\infty} \frac{X_n(nt)}n
\le \limsup_{n\to\infty} \frac{X_n(nt)}n \le x^+(b,t)
\tag 2.21
$$
almost surely. In particular, if $x(b,t)=x^\pm(b,t)$,
then we have a strong law of large numbers: 
$$
\lim_{n\to\infty} \frac{X_n(nt)}n = x(b,t)\quad\text{a.s.} 
\tag 2.22
$$
\endproclaim

In the special case where characteristics are unique
we can dispense with assumption  (2.20). 
The second class particle tracks the 
macroscopically defined characteristic that emanates 
from its  random initial position, as stated in part (a)
of the next theorem. We only need to
assume that  the fluctuations of the initial location
are macroscopically bounded as $n$ grows:
$$
-\infty< \liminf_{n\to\infty} \frac{X_n(0)}n
\le \limsup_{n\to\infty} \frac{X_n(0)}n <\infty
\qquad\text{almost surely.} 
\tag 2.23
$$
The uniqueness of characteristics is true if we look at
characteristics started at a time $s>0$. But only
if the flux function is differentiable, 
so this assumption needs to be made explicitly for 
 $2\le K<\infty$ in part (b).

\proclaim{Theorem 4} Assume {\rm (2.3)} and  {\rm (2.23)}. 

{\rm (a)} If 
$x(y,t)=x^\pm(y,t)$ for all $y\in\bold R$, then 
almost surely
$$
\lim_{n\to\infty} \left| \frac{X_n(nt)}n
- x\left( \frac{X_n(0)}n \,,\,t\right) \right|=0.
\tag 2.24
$$

{\rm (b)} Let $0<s<t$. 
 In the case $2\le K<\infty$
suppose $f$ is differentiable
on $[0,K]$. For 
  $K=1$ and $K=\infty$ this is known and does not need
to be assumed.
Then $x(y;s,t)=x^\pm(y;s,t)$ for all $y\in\bold R$,
and almost surely
$$
\lim_{n\to\infty} \left| \frac{X_n(nt)}n
- x\left( \frac{X_n(ns)}n \,;s,\,t\right) \right|=0. 
\tag 2.25
$$
\endproclaim

The meaning of $x\left( n^{-1}{X_n(0)} \,,\,t\right)$ 
is that the deterministic function $x(\cdot,t)$ defined 
by (2.18--19) is evaluated at the random point $n^{-1}{X_n(0)}$.
If assumption (2.23) is strengthened to  
$na\le X_n(0)\le nb$ for all large enough $n$, 
then for Theorem 4(a) it is sufficient to assume
$x(y,t)=x^\pm(y,t)$ for $a\le y \le b$. As for 
Theorem 1, if the assumptions are valid 
in probability, then so are the conclusions. 

\hbox{}

\subhead 2.3. Comments on the assumptions and results
\endsubhead

\hbox{}

{\it 2.3.1. Differentiability and strict convexity.}
Since $g$ is convex,  left and right derivatives 
$$g'_-(x)=\lim_{\e\searrow 0}\frac{g(x-\e)-g(x)}{-\e}
\qquad\text{and}\qquad
g'_+(x)=\lim_{\e\searrow 0}\frac{g(x+\e)-g(x)}{\e}
%\tag 
$$
exist at all points $x$. Same is true for $f$, and 
for both $f$ and $g$ differentiability on an open 
interval is equivalent to continuous differentiability. 
 {\it Strict convexity} means that
there are no $x_1<x_2$ such that $g'_+(x_1)=g'_-(x_2)$,
which is the same as not having a line segment in the
graph of $g$. 
From the duality (2.8) it follows that $f$ is differentiable
on $[0,K]$ iff $g$ is strictly convex on $[-1,1]$, 
and $f$ is strictly concave on $[0,K]$ iff $g$ is 
differentiable on $\bold R$. Consult Rockafellar (1970) 
for a general treatment of convex analysis. 

 Theorem 4(b) is true because strict convexity 
of $g$ guarantees that 
$x(y;s,t)$ $=$ $x^\pm(y;s,t)$ for all $y\in\bold R$, if $s>0$.
This is proved in Proposition 3.1 in  Section 3.2.

The restriction $x-(t-s)\le y\le x+(t-s)$
 in definition (2.16--17) of
$y^\pm(x;s,t)$ is necessary because $g$ is linear
outside $[-1,1]$.  For example, 
take as initial profile 
$\rho_0(x)=K$ for $x\le 0$ and $\rho_0(x)=0$ for $x> 0$. 
Then for $x>t$, the entire interval $[0,x-t]$ 
consists of maximizers in the Hopf-Lax formula. 
Without the restriction in (2.16) we would 
have $y^-(x,t)=0$ for all $x>t$, and then  
$x^+(0,t)=\infty$ which is not the appropriate characteristic. 

\hbox{}

{\it 2.3.2. Shocks.} Consider the simplest shock case
(Riemann problem)
with initial profile 
$\rho_0(x)=\alpha\bold 1\{x<0\}+\beta\bold 1\{x>0\}$
where $\alpha<\beta$. The solution (2.12) is 
$\rho(x,t)=\alpha\bold 1\{x<\xi t\}+\beta\bold 1\{x>\xi t\}$
with shock speed $\xi=(f(\beta)-f(\alpha))/(\beta-\alpha)$. 
The unique characteristic from the origin is 
$x(0,t)=\xi t$.  Theorem 3 says 
that the second class particle converges to the 
shock. This generalizes results of Ferrari (1992) 
in two ways: From exclusion to  $K$-exclusion, and 
from initial product distributions to general
initial distributions that have the shock  
profile.

\hbox{}

{\it 2.3.3. Rarefaction fan.} 
Ferrari and Kipnis  (1995) proved 
convergence in probability in (2.25) for the following
case of TASEP: initial distribution is 
 product measure with density $\rho$ to the
left and density $\lambda$ to the right of the origin, 
with $\rho>\lambda$, and the second class particle 
initially at the origin. This initial profile 
$\rho_0(x)=\rho\bold 1\{x<0\}+\lambda\bold 1\{x>0\}$
has an inadmissible shock at the origin, which
produces a ``rarefaction fan'' of characteristics 
with $x^-(0,t)=(1-2\rho)t$ and 
$x^+(0,t)=(1-2\lambda)t$. Ferrari and Kipnis also
proved that weakly $n^{-1}X_n(nt)$ converges to 
the uniform distribution on
$[x^-(0,t),x^+(0,t)]$. Consequently one cannot hope
for a law of large numbers when 
$x^-(0,t)<x^+(0,t)$. 

If $f$ is not differentiable at $\rho$,
then 
the situation $x^-(b,t)<x^+(b,t)$ happens even for the
constant profile $\rho(x,t)\equiv \rho$. 
See Example 3.1 in Section 3.2. In this case 
we cannot deduce convergence of $n^{-1}X_n(nt)$ 
from
Theorem 3, only the bounds (2.24), unless we can prove
the differentiability of $f$.

\hbox{}

{\it 2.3.4. The characteristic ordinary differential equation.}
When $f$ is strictly concave
and  differentiable, Theorem 3 also contains
this statement: The characteristics 
$x(t)=x^\pm(b,t)$ are {\it Filippov solutions} of the o.d.e.
$$
\frac{dx}{dt} =f'(\rho(x,t))\;,\;x(0)=b\,,
\tag 2.26
$$
where $\rho(x,t)$ is the Lax-Oleinik solution of (2.5).
Since $\rho(x,t)$ does not exist at all $x$, (2.26) is 
interpreted as 
$$\essinf_x f'(\rho(x,t)) \le \frac{dx}{dt}
\le \esssup_x f'(\rho(x,t))\qquad\text{for a.e.\ $t$.}
\tag 2.27
$$
The essential infimum relative
to $x$, at the point $(x,t)$, is defined by
$$\essinf_x f'(\rho(x,t))
=\lim_{\delta\searrow 0}\; \essinf_{x-\delta<y<x+\delta}
f'(\rho(y,t)),
$$
and similarly for the essential supremum.  The 
Lax-Oleinik formula (we review it in Section 3.1 below)
shows that (2.27) is equivalent to 
$$f'\left(-g'\left(\frac{x-y^+(x,t)}t\right)\right)
\le \frac{dx}{dt} \le 
f'\left(-g'\left(\frac{x-y^-(x,t)}t\right)\right)
\tag 2.28
$$
for almost every $t$. 
These results are developed  in Rezakhanlou (1995) for 
the exclusion and zero-range processes. 

At present 
 we cannot ascertain strict concavity or
 differentiability of $f$ 
for $2\le K<\infty$. This is a problem,
because to interpret (2.26) in the Filippov sense, 
 $f'(\rho(x,t))$ must exist at least a.e. But
if $f'$ does not exist at $\rho$ and we consider 
the constant solution $\rho(x,t)\equiv \rho$, then 
 $f'(\rho(x,t))$ does not exist at any $(x,t)$.  

Without any assumptions on $f$ beyond concavity, 
 bounds (2.28) are valid
for  $x(t)=x^\pm(b,t)$
 in the form 
$$f_+'\left(-g_+'\left(\frac{x-y^+(x,t)}t\right)\right)
\le \frac{dx}{dt} \le 
f_-'\left(-g_-'\left(\frac{x-y^-(x,t)}t\right)\right)
\tag 2.29
$$
for almost every $t$. If we assume $f$ differentiable, 
then in fact (2.26)
 is valid again in the sense (2.27), with $\rho(x,t)$
given by (2.12). In other words, the o.d.e.\ works
 even though the 
Lax-Oleinik formula fails without strict concavity of $f$. 
We omit the proof  of (2.29) and the last statement
above.

\hbox{}

\subhead 2.4. Generalization: particles jump in batches
\endsubhead
We discuss here a generalization of the 
$K$-exclusion and zero-range processes for which the
 variational coupling technique works and Theorems 1--4
 are true. 

From a queueing perspective, the particles are
customers moving through  an infinite
sequence of servers. Each server has space for $K$ customers
in its queue, and $\eta(i,t)$ is the number
of customers waiting at server $i$ at time $t$.
In the formulation (2.1), each server serves at exponential
rate 1, and only one customer is served at a time
and moves on to the next server. 
We generalize this by permitting customers to move 
 in batches of varying sizes, with exponential 
rates and independently for each batch size.  
Let $\beta_h$ be the rate at which a batch of 
$h$ customers can jump from $i$ to $i+1$, 
$h=1,2,3,\ldots$ . 
If the  state just before the jump 
 does not permit $h$ customers
to jump [so either $\eta(i)<h$ or $\eta(i+1)>K-h$], 
we do not suppress the entire jump but instead move as 
many customers as possible without violating the 
constraints. 

The family of processes is parametrized by
 $K$ and 
 a sequence of nonnegative numbers $(\beta_h: h\in\bold N)$. 
If $K<\infty$ assume $\beta_h=0$ for $h>K$.
To realize this dynamics take a collection of mutually independent 
Poisson point processes $\{D^h_i:i\in\bold Z, h\in\bold N\}$
on the time axis
$(0,\infty)$, so that the rate of $D^h_i$ is $\beta_h$. 
At each epoch of $D^h_i$, a batch of $h$ particles
(customers) attempt to jump
from 
site $i$ to site $i+1$. 
The generator becomes 
$$Lf(\eta)=\sum_{i\in\bold Z} \sum_{h\in\bold N} \beta_h
[f(\eta^{h,i,i+1})-f(\eta)],
%\tag 2.1
$$
where $\eta^{h,i,i+1}$ is the configuration that results from an
attempt to move a batch
of $h$ particles  from site $i$ to $i+1$. Denote
the actual 
number of particles that move by
$$b=b(h,\eta(i),\eta(i+1))=\min\{h,\, \eta(i), \,K-\eta(i+1)\}. 
%\tag 
$$
Then the new configuration is 
$$ 
\eta^{h,i,i+1}(j)=
\cases \eta(i)-b, &j=i\\
\eta(i+1)+b, &j=i+1\\
\eta(j), &j\ne i,i+1.
\endcases
$$

The definition of a second class particle is the same 
as before. From the queueing perspective 
it is  a second class customer. He
can move forward only if first class customers
 do not fill the batch that is served,
and he has to move back to the 
previous server if arriving first class customers fill the 
capacity of the queue. 

The hydrodynamic limit Theorem 1 is valid as stated
for the batch process too, provided the jump rates
are restricted so that we can expect finite limits.
It is enough to require 
$\sum_{h\in\bold N} \beta_hh^2<\infty $.
Under this assumption a weaker version of Theorem 1 
with convergence in probability was proved in
Sepp\"al\"ainen (2000). That paper also shows how
to adapt the variational 
coupling approach to batch jumps. 

 If  $K=\infty$, the i.i.d.\ geometric probability measures
 (2.2) 
 are again  invariant, for any choice of the rates $(\beta_h)$.
In this case the flux function is 
$$f(\rho)=\sum_{h\in\bold N}
\beta_h\left[ \frac{\rho}{1+\rho}
+ \left(\frac{\rho}{1+\rho} \right) +\cdots +
\left(\frac{\rho}{1+\rho} \right)^h
\right].
%\tag 
$$

\hbox{}

\head 3. Scalar conservation laws with concave flux functions  \endhead 

In this section we investigate the consequences
of the lack of  strict concavity and 
differentiability of the flux function. It is important
to clarify these facts because they directly influence 
the results we obtain for the particle models. 
Section 3.1 looks at strict concavity and  section 3.2 
at differentiability. In section 3.3 we prove the 
uniqueness criterion of Theorem 2 which assumes neither
strict concavity nor differentiability of the flux.

\subhead 3.1. Strict concavity of the flux and entropy
shocks   \endsubhead
We begin by recalling the Lax-Oleinik formula
from  Section 3.4 in Evans (1998),
and then show how it depends on the strict concavity 
of the flux function. 

\demo{Lax-Oleinik formula}
Suppose $f$ is strictly concave and 
differentiable on its domain $[0,K]$.  Define two functions 
$$
\rho^\pm(x,t)=-g'\left(\frac{x-y^\pm(x,t)}t\right).
\tag 3.1
$$
For fixed $t>0$, under these assumptions $y^+(x,t)=y^-(x,t)$
for all but countably many $x$, and consequently 
the Lax-Oleinik solution  
$$\rho_{\lo}(x,t)=\rho^\pm(x,t)
\tag 3.2
$$ is well-defined 
for all but countably many $x$. 
It satisfies $\rho_{\lo}(x,t)=u_x(x,t)$
 a.e., and is a weak solution of (2.5) in the integral sense (2.13). 
\enddemo

In particular, when it exists
$\rho_{\lo}$ 
is the solution  defined in (2.12) for the hydrodynamic limit. 
One can check that  $\rho_{\lo}(\cdot\,,t)$ is continuous
at all but the countably many $x$ at which 
$\rho^-(x,t)\ne \rho^+(x,t)$. 
If 
 $\rho_{\lo}(\cdot, t)$ is discontinuous at $x$, 
there must be a jump upward:  
$\rho_{\lo}(x-,t)< \rho_{\lo}(x+,t). 
$
This is because  $g'$ is nondecreasing and
$y^-(x,t)\le y^+(x,t)$. 
 Let us call such a jump
an {\it entropy shock} for a concave flux function. 
[For a convex flux entropy shocks are downward jumps.]

There are two ways definition (3.1)--(3.2) 
could conceivably run into a problem: 

(i) If $f$ 
is not strictly concave then $g'$ fails to exist at some
points and the right-hand side 
of (3.1) may not be defined.

 (ii) If $f'$ does not exist everywhere then 
$y^-(x,t)< y^+(x,t)$ may happen for a nontrivial 
interval of points $x$, and formula (3.2) is not true a.e. 
See Example 3.1 in Section 3.2.

It turns out that  if
 $f$ is strictly concave, 
definition (3.1--2) works even if  $y^-(x,t)= y^+(x,t)$
fails. And all discontinuities are entropy shocks.

\proclaim{Theorem 5} Suppose the flux $f$ is strictly
concave. Then for a fixed $t>0$ 
 the solution $\rho(x,t)=u_x(x,t)$ 
defined as in {\rm (2.12)} exists and is
continuous  on a set $H\subseteq\bold R$
that contains all but countably many $x$. On
$H$ $\rho(x,t)=\rho^\pm(x,t)$ so   the Lax-Oleinik formula
is valid.  For all $x$ 
$\rho(x-,t)\le \rho(x+,t)$, so 
the discontinuities of $\rho(\cdot, t)$ are entropy
shocks. 

Conversely, suppose $f$ has a linear segment.
Then there is a 
range of densities $\rho$ such that the constant solution
$\rho(x,t)\equiv\rho$ cannot be represented by 
any left or right derivative of
$g$, so the Lax-Oleinik formula {\rm (3.1--2)} fails. Furthermore, 
the relevant solution defined by {\rm (2.12)} can have 
 nonentropy shocks $\rho(x-,t)> \rho(x+,t)$. 
\endproclaim

\demo{Proof} 
 Strict concavity of 
 $f$ implies  that $g$ is differentiable. 
Inequalities (3.3)--(3.4) below then imply that the right and left
$x$-derivatives of $u(x,t)$ exist and are given by
$$u_{x\pm}(x,t)=-g'\left(\frac{x-y^\pm(x,t)}t\right)=\rho^\pm(x,t).$$
Since $g'$ is continuous and $y^+(\cdot\,,t)$ 
nondecreasing and right-continuous, we
see that $\rho^+(\cdot\,,t)$ is right-continuous, has at most 
countably many discontinuities,  and has left limits at all
points. And correspondingly for $\rho^-(x,t)$.  

By the Lipschitz property $u(x,t)$ is differentiable at 
a.e.\ $x$ and is the integral of its derivative. Thus
$$u(x,t)=\int_0^x u_x(\xi,t)d\xi=\int_0^x \rho^\pm(\xi,t)d\xi.$$
Consequently $\rho(x,t)\equiv u_x(x,t)$ exists and equals 
$\rho^\pm(x,t)$ at any $x$ where  one of $\rho^\pm(x,t)$
 is continuous. 
This defines the set $H$. From this follows 
$\rho(x-,t)=\rho^-(x,t)\le \rho^+(x,t)=\rho(x+,t)$. 

For the converse, suppose $f'(\rho)=\xi$ for $\rho_0<\rho<\rho_1$.
Then $g$ has a corner at $\xi$, with 
$g'_-(\xi)=-\rho_1$ and $g'_+(\xi)=-\rho_0$. 
The values $\rho\in(\rho_0,\rho_1)$ are not taken
by $g'_\pm$ at any point, so a 
constant solution  $\rho(x,t)\equiv\rho\in(\rho_0,\rho_1)$ 
 cannot be 
represented 
 by a derivative of $g$. 

The initial profile 
$\rho_0(x)=K$ for $x\le 0$ and $\rho_0(x)=0$ for $x> 0$
gives the solution $u(x,t)=-tg(x/t)$. Then 
$\rho((t\xi)-,t)=\rho_1>\rho_0=\rho((t\xi)+,t)$ is a 
nonentropy shock. 
\qed
\enddemo

To finish this section we prove the inequalities used 
in the proof above. 
Let us write $\displaystyle{\overline{\varliminf_{\e\searrow 0}}}$
to simultaneously include both the $\limsup$ and the $\liminf$
as $\e\searrow 0$.

\proclaim{Lemma 3.1} The following inequalities hold for 
right and left $x$-derivatives of $u(x,t)$ at fixed $t>0$. 
$$
-g'_+\left(\frac{x-y^-(x,t)}{t}\right)
\le 
\overline{\varliminf_{\e\searrow 0}}\;
\frac{u(x-\e,t)-u(x,t)}{-\e}
\le 
-g'_-\left(\frac{x-y^-(x,t)}{t}\right)
\tag 3.3
$$
and 
$$
-g'_+\left(\frac{x-y^+(x,t)}{t}\right)
\le 
\overline{\varliminf_{\e\searrow 0}}\;
\frac{u(x+\e,t)-u(x,t)}{\e}
\le 
-g'_-\left(\frac{x-y^+(x,t)}{t}\right). 
\tag 3.4
$$
\endproclaim

\demo{Proof} The proofs use convexity, the Hopf-Lax 
formula, and the right [left] continuity
of $y^+(\cdot\,,t)$ and $g'_+$ [$y^-(\cdot\,,t)$ and $g'_-$]. 
As an illustration, here is the argument for 
the second inequality of (3.4). We leave the rest to the reader. 
$$\aligned
&u(x+\e,t)-u(x,t)\\
&\le u_0(y^+(x+\e,t)) -tg\left(\frac{x+\e-y^+(x+\e,t)}{t}\right)\\
&\qquad\qquad -u_0(y^+(x+\e,t)) +tg\left(\frac{x-y^+(x+\e,t)}{t}\right)\\
&= -\e\cdot \left(\frac{t}{\e}\right)
\cdot \left\{g\left(\frac{x-y^+(x+\e,t)}{t}+\frac{\e}{t}\right) 
-g\left(\frac{x-y^+(x+\e,t)}{t}\right) \right\}\\
&\le -\e\cdot g'_+\left(\frac{x-y^+(x+\e,t)}{t}\right).  
\endaligned
$$
Consider two possibilities: If $y^+(x+\e,t)=y^+(x,t)$ for 
small enough $\e>0$, then 
$$ g'_+\left(\frac{x-y^+(x+\e,t)}{t}\right)
=g'_+\left(\frac{x-y^+(x,t)}{t}\right) \ge
g'_-\left(\frac{x-y^+(x,t)}{t}\right).
$$
The other possibility is that $y^+(x+\e,t)>y^+(x,t)$ for 
all $\e>0$. From right continuity of $x\mapsto y^+(x,t)$
it follows that 
$(x-y^+(x+\e,t))/{t}$ increases strictly to $(x-y^+(x,t))/{t}$
as $\e\searrow 0$. And then 
$$\lim_{\e\searrow 0} g'_+\left(\frac{x-y^+(x+\e,t)}{t}\right)
= g'_-\left(\frac{x-y^+(x,t)}{t}\right).
$$
In either case we obtain  the second inequality of (3.4). 
\qed
\enddemo

\subhead 3.2. Differentiability of the flux and 
characteristics
\endsubhead
In Corollary 3.1(b) below 
 we show that differentiability of $f$ implies
$y^-(x,t)=y^+(x,t)$ for all but countably many $x$, 
the property used in part (3.2) of the Lax-Oleinik formula.
In Proposition 3.1
 we prove the property needed for Theorem 4(b), namely
that $x^-(b;s,t)=x^+(b;s,t)$ for $0<s<t$. 
Example 3.1 shows how various properties can fail
without differentiability of $f$. 
The hypothesis of differentiability of $f$ is used
 in the equivalent form 
that $g$ is strictly convex on $[-1,1]$. 

Let $I(x;s,t)$ denote the set of $y\in[x-(t-s),x+(t-s)]$ at which
the supremum in (2.15) is achieved. 

\proclaim{Lemma 3.2} Suppose $0\le t_1<t_2<t_3$, $y_2\in I(x;t_2,t_3)$,
and $y_1\in I(y_2;t_1,t_2)$. Then $y_1\in I(x;t_1,t_3)$. 
And either 
$$\frac{x-y_1}{t_3-t_1}=\frac{x-y_2}{t_3-t_2}=\frac{y_2-y_1}{t_2-t_1}
$$
or $g$ is linear between 
$(y_2-y_1)/(t_2-t_1)$ and $(x-y_2)/(t_3-t_2)$.
\endproclaim

\demo{Proof} Modify the proof of 
Rezakhanlou's (1995) Lemma 3.1. 
 \qed
\enddemo

\proclaim{Lemma 3.3} Suppose $g$ is strictly convex
on $[-1,1]$, 
$x_1<x_2$, $y_i\in I(x_i;s,t)$ for $i=1,2$. Then $y_1\le y_2$. 
\endproclaim

\demo{Proof}  The assumptions  $y_i\in I(x_i;s,t)$ give 
$$
u(y_1,s)-(t-s)g\biggl(\frac{x_1-y_1}{t-s}\biggr)
\ge 
u(y_2,s)-(t-s)g\biggl(\frac{x_1-y_2}{t-s}\biggr)
$$
and 
$$
u(y_2,s)-(t-s)g\biggl(\frac{x_2-y_2}{t-s}\biggr)
\ge 
u(y_1,s)-(t-s)g\biggl(\frac{x_2-y_1}{t-s}\biggr). 
$$
To get a contradiction, suppose
$y_1-y_2>0$. Then the above inequalities give 
$$
\frac{g\left(\frac{x_2-y_2}{t-s}\right)
-g\left(\frac{x_2-y_1}{t-s}\right)}
{\left(\frac{x_2-y_2}{t-s}\right)
-\left(\frac{x_2-y_1}{t-s}\right)}
\le
\frac{g\left(\frac{x_1-y_2}{t-s}\right)
-g\left(\frac{x_1-y_1}{t-s}\right)}
{\left(\frac{x_1-y_2}{t-s}\right)
-\left(\frac{x_1-y_1}{t-s}\right)}\,.
$$
Since $x_2-y_2>x_1-y_2$,  $x_2-y_1>x_1-y_1$,  and
$(x_i-y_i)/(t-s)\in[-1,1]$, 
this contradicts the assumption on $g$. 
\qed
\enddemo

\proclaim{Corollary 3.1} Suppose $g$ is strictly convex on $[-1,1]$.

{\rm (a)}  For $x_1<x_2$, $y^+(x_1;s,t)\le y^-(x_2;s,t)$. 

{\rm (b)}  For fixed $0\le s<t$, $y^-(x;s,t)= y^+(x;s,t)$
 for
 all but countably many $x$. 

{\rm (c)}  $x^+(b;s,t)=\inf\{ x: y^+(x;s,t)> b\}$ and 
$x^-(b;s,t)=\sup\{ x: y^-(x;s,t)< b\}$. 

{\rm (d)} If $x^-(b;s,t)<z<x^+(b;s,t)$, then $y^\pm(z;s,t)=b$. 
\endproclaim

\demo{Proof} (a) is immediate from Lemma 3.3. 

(b) Suppose
$x_0$ is a continuity point of the nondecreasing, 
right-continuous function $x\mapsto y^+(x;s,t)$. 
Then  
$$y^+(x_0;s,t)\ge y^-(x_0;s,t)\ge \lim_{x\nearrow x_0}
y^+(x;s,t)= y^+(x_0;s,t), 
$$
where the second inequality follows from (a). 
Thus $y^-(x;s,t)= y^+(x;s,t)$ at all continuity points
of $y^+(x;s,t)$, which includes all but countably many $x$. 

(c) Let $w(b)=\inf\{ x: y^+(x;s,t)> b\}$. We have 
$x^+(b;s,t)\ge w(b)$ without any new assumptions: 
if $x<w(b)$, then 
$$y^+(x;s,t)\le b \ \Longrightarrow\  y^-(x;s,t)\le b
\ \Longrightarrow\  x^+(b;s,t)\ge x.
$$
For the converse we need part (a). Let $x>w(b)$. 
Then there exists $x_1$ such that $x>x_1>w(b)$. 
$$y^+(x_1;s,t)> b\  \Longrightarrow\  y^-(x;s,t)> b \ \Longrightarrow\  
 x^+(b;s,t)\le x.
$$
We leave the proof of the other formula to the reader.

(d) Immediate from definition (2.18--19) and part (c). 
\qed
\enddemo

\proclaim{Proposition 3.1} Suppose $g$ is strictly convex
on $[-1,1]$,
and $s>0$. Then $x^-(b;s,t)=x^+(b;s,t)$ for all $b\in\bold R$
and $t>s$. 
\endproclaim

\demo{Proof} Suppose there exists a $z$ such
that $x^-(b;s,t)<z<x^+(b;s,t)$. By Corollary (d) 
$y^\pm(z;s,t)=b$, and consequently  for $\e>0$, 
$$u(b,s)-(t-s)g\left(\frac{z-b}{t-s}\right)
\ge 
u(b\pm\e,s)-(t-s)g\left(\frac{z-b\mp \e}{t-s}\right).
$$
Letting $\e\searrow 0$ then gives
$$\limsup_{\e\searrow 0}\frac{u(b+\e,s)-u(b,s)}{\e}
\le 
-g'_-\left(\frac{z-b}{t-s}\right)
$$
and 
$$\liminf_{\e\searrow 0}\frac{u(b-\e,s)-u(b,s)}{-\e}
\ge 
-g'_+\left(\frac{z-b}{t-s}\right).
$$
Since  we may vary $z$ in an interval, and since 
by assumption the slope of $g$ cannot remain
constant in this interval, we conclude the strict
inequality 
$$\limsup_{\e\searrow 0}\frac{u(b+\e,s)-u(b,s)}{\e}
<
\liminf_{\e\searrow 0}\frac{u(b-\e,s)-u(b,s)}{-\e}.
$$
By (3.3--4), this implies
$$g'_+\left(\frac{b-y^+(b,s)}{s}\right)
>
g'_-\left(\frac{b-y^-(b,s)}{s}\right).
\tag 3.5
$$
The slope of a convex function is nondecreasing
and by definition $(b-y^-(b,s))/{s}\ge (b-y^+(b,s))/{s}$. 
So  (3.5) implies that $y^-(b,s)=y^+(b,s)$ 
and $g$ has a corner at $(b-y^\pm(b,s))/{s}$. 
Since we are assuming $g$ strictly convex  in $[-1,1]$, Lemma 3.2
now forces $(z-b)/(t-s)=(b-y^\pm(b,s))/s$. This cannot
hold for distinct $z$'s. But 
 the  preceding argument is valid for all 
$z\in(x^-(b;s,t),x^+(b;s,t))$, so to avoid
contradiction  
we must have $x^-(b;s,t)=x^+(b;s,t)$.
\qed
\enddemo

Finally we show by example how 
a point of nondifferentiability in $f$ destroys all 
the good properties of the characteristics. 

\demo{Example 3.1} 
Suppose $\rho\in(0,K)$ is a point such that
 $f'_-(\rho)=\xi_2>\xi_1=f'_+(\rho)$. Then 
$g$ has a linear segment of slope $-\rho$ 
over $[\xi_1,\xi_2]$. Apply the formulas 
 to the constant initial profile $\rho_0(x)\equiv \rho$,
with $u_0(x)=\rho x$. We get 
$I(x;s,t)=[x-\xi_2(t-s),x-\xi_1(t-s)]$, 
$y^-(x;s,t)=x-\xi_2(t-s)$, $y^+(x;s,t)=x-\xi_1(t-s)$,
$x^-(b;s,t)=b+\xi_1(t-s)$, $x^+(b;s,t)=b+\xi_2(t-s)$,
and 
$u(x,t)=\rho x-tf(\rho)$.  In particular 
$\rho(x,t)=u_x(x,t)\equiv \rho$, so the constant
profile is the solution. But the 
conclusions of  Lemma 3.3,  
Corollary 3.1(a)--(b), and Proposition 3.1 have been contradicted. 
\enddemo

\subhead 3.3. Proof of Theorem 2 \endsubhead
We prove the uniqueness criterion of Theorem 2
under these  general assumptions: The 
flux $f$  is continuous and 
concave on an interval $I$, and $g$ is defined
by $g(x)=\sup_{\rho\in I}\{f(\rho)-\rho x\}$.
  The initial density
 $\rho_0(x)$ is a given locally bounded measurable
function on $\bold R$,   $u(x,t)$ is a locally Lipschitz 
function on $\bold R\times[0,\infty)$
 defined by the Hopf-Lax formula (2.10), and 
$\rho(x,t)=u_x(x,t)$ exists for a.e. $x$ for any fixed $t$.  
$\lambda(x,t)$ is a locally bounded measurable function on 
$\bold R\times [0,\infty)$ that  satisfies the 
boundary condition $\lambda(x,0)=\rho_0(x)$ a.e.\ and 
the integral condition (2.13). The
 ranges of $\rho$ and $\lambda$ lie in $I$. 

 The goal is to show that, for a.e.\ $x$, 
$$\int_0^t f(\lambda(x,s))ds\le \int_0^t f(\rho(x,s))ds  
\tag 3.6
$$
 for all $t\ge 0$. 
And conversely, if 
 equality holds in (3.6) for
almost all $(x,t)$, then $\lambda(x,t)=\rho(x,t)$ 
for almost all $(x,t)$.

Some technical difficulties stem from 
 assuming only (2.13) and no regularity on $\lambda$.
We cannot expect the function $x\mapsto \lambda(x,t)$ to be a 
sensible object for every fixed $t$
because $\lambda(\cdot,t)$ can be redefined on any
 null set of $t$'s without affecting (2.13). 
 At $t=0$ we have 
$$
\lim_{\e\to 0}\frac1\e \int_0^\e dt \int_{\bold R}
\lambda(x,t)\psi(x)dx = \int_{\bold R} \rho_0(x)\psi(x)dx
\tag 3.7
$$
for $\psi\in C_c^1(\bold R)$.
This is obtained from (2.13) by taking $\phi(x,t)=\psi(x)g(t)$
where $g(0)=1$, $g'=-1/\e$ on $(0,\e)$,
and $g(t)=0$ for $t\ge \e$.  [The $C^1$ test functions in (2.13)
can be replaced by continuous, piecewise $C^1$ test 
functions by taking limits.]

A corresponding property holds for certain later times $t$. 
 Abbreviate
$$\leb(h)=\text{the set of Lebesgue points of the function $h$}
$$
for any measurable, locally bounded function $h$ on either  
$\bold R$ or $\bold R\times[0,\infty)$. [For the definition 
and properties of Lebesgue points see Section 7.2 in
Wheeden and Zygmund (1977).] Define 
$$
\bold T=\{t>0: \text{$(x,t)\in\leb(\lambda)$
for almost every $x\in\bold R$}\}.
$$
Since the complement of $\leb(\lambda)$ is a null set
of $\bold R\times(0,\infty)$, almost every $t$ lies in $\bold T$. 
By considering test functions of the type $\phi(x,t)=\psi(x)g(t)$
for suitable $g$, one obtains from (2.13) a limit such as (3.7) around
$t\in\bold T$, and then the following for $t\in\bold T$:
$$\int_{\bold R} \psi(x)\lambda(x,t)dx
-\int_{\bold R} \psi(x)\rho_0(x)dx
=\int_0^t \int_{\bold R} f(\lambda(x,s))\psi'(x)dx\, ds
\tag 3.8
$$
for compactly supported $C^1$ functions $\psi$ on $\bold R$.
By taking limits of such functions, (3.8) is valid
 for 
compactly supported, continuous, piecewise $C^1$ functions $\psi$.

\proclaim{Lemma 3.4} 
{\rm (i)} Fix finite $A<B$ and $T>0$. Then there exists a 
constant $C=C(A,B,T)$ such that,   
for all  $A\le a<b\le B$ and $s,t\in\bold T\cup\{0\}$
such that $0\le s,t\le T$, 
$$\left| \int_a^b \lambda(x,t)dx-\int_a^b \lambda(x,s)dx\right|
\le C\cdot |t-s|.
\tag 3.9
$$

{\rm (ii)} If $t\in\bold T\cup\{0\}$ and
 $(a,t),(b,t)\in\leb(f\circ\lambda)$, then we have 
a time derivative along points of $\bold T$:
$$\lim\Sb s\to t\\s\in\bold T\endSb \frac1{t-s}\left\{
\int_a^b \lambda(x,t)dx-\int_a^b \lambda(x,s)dx\right\}
=f(\lambda(a,t))-f(\lambda(b,t)).
\tag 3.10
$$
\endproclaim

\demo{Proof} In (3.8) take $\psi$ such that 
$\psi=0$  outside $[a-\e,b+\e]$, $\psi'=1/\e$
on 
$(a-\e,a)$, 
 $\psi=1$ on $[a,b]$, and $\psi'=-1/\e$ on 
 $(b,b+\e)$. We get
$$\aligned
& \int_a^b \lambda(x,t)dx-\int_a^b \lambda(x,s)dx
+O(\e)\\
=\;&
\frac1\e\int_s^t \int_{a-\e}^a f(\lambda(x,\tau))dxd\tau
-\frac1\e\int_s^t \int_b^{b+\e}f(\lambda(x,\tau))dxd\tau.
\endaligned
\tag 3.11
$$
The $O(\e)$ term contains the integrals over $[a-\e,a]$ and 
$[b,b+\e]$ from the left-hand side. These are bounded by $O(\e)$ 
uniformly  over the range of $a,b,s,t$ by the local 
boundedness of $\lambda$.   The $1/\e$ factors
on the right are the slopes of $\psi$. 
 The second line in (3.11)
is bounded by $ C |t-s|$ in absolute value, where 
$C$ is twice the supremum of $|f(\lambda(x,\tau))|$
 over $[A,B]\times[0,T]$. 
Part (i) follows by letting $\e\searrow 0$. For
part (ii) set $\e=\delta|t-s|$, divide through (3.11) by $t-s$,
 let first $\bold T\ni s\to t$ 
and then $\delta\searrow 0$. 
\qed
\enddemo

As a locally Lipschitz function $u(x,t)$ is differentiable 
a.e.\  by Rademacher's theorem [Section 5.8.3 in
Evans (1998)]. And as the viscosity solution defined 
by the Hopf-Lax formula, $u(x,t)$ satisfies the equation
$u_t+f(u_x)=0$ a.e.\ [Theorem 5, Section 3.3.2
in Evans (1998)]. By definition 
$\rho=u_x$ exists a.e., so  for almost every $x_0$
and all $t>0$, 
$$\int_0^t f(\rho(x_0,s))ds=-\int_0^t u_t(x_0,s)ds
=u_0(x_0)-u(x_0,t).
\tag 3.12
$$

It is also the case that for almost every $x_0$, 
$(x_0,t)\in\leb(f\circ\lambda)$ for almost
every $t$. 
Fix such an $x_0$ for which also (3.12) holds.  
 For $(x,t)\in\bold R\times \bold T$ define
$$v(x,t)=\int_{x_0}^x \lambda(y,t)dy -\int_0^t f(\lambda(x_0,s))ds
+u_0(x_0).
\tag 3.13
$$
The first integral is interpreted with a sign so that 
$v(x_2,t)-v(x_1,t)=\int_{x_1}^{x_2} \lambda(y,t)dy$ for 
$x_1<x_2$.  By Lemma 3.4(i) and local
boundedness of $\lambda$ and $f\circ\lambda$, 
 $v$ is locally Lipschitz on $\bold R\times \bold T$.
Consequently $v$ extends uniquely to a 
locally Lipschitz function on $\bold R\times [0,\infty)$.
This extension has the correct boundary 
values at $t=0$ because by Lemma 3.4(i), 
$$\aligned
\lim\Sb \bold T\ni t\to 0\endSb v(x,t)
&=\int_{x_0}^x \lambda(y,0)dy +u_0(x_0)
=\int_{x_0}^x \rho_0(y)dy +u_0(x_0)\\
&=u_0(x).
\endaligned
$$

The goal is now to prove 
$$\text{$v(x,t)\ge u(x,t)$ for all $(x,t)$.}
\tag 3.14
$$
Setting $x=x_0$ in (3.14) then gives, by (3.12) and (3.13),
  the conclusion (3.6)
 for $x=x_0$ and all $t\in\bold T$. Both sides of (3.6) are
continuous in $t$ because the integrands are locally bounded,
so (3.6) follows for all $t$. This whole argument beginning
with (3.12) can be repeated for almost every $x_0$. 
At the end of this section we argue
that equality in (3.6) implies $\lambda=\rho$ a.e. 
First we prove 
(3.14).

Fix $(x,t)$. For a fixed $y\in\bold R$ set 
$\xi(s)=y+(s/t)(x-y)$ for $s\in[0,t]$. As $y$ varies
we consider the 
collection of line segments  $\{(\xi(s),s): 0\le s\le t\}$
on the plane $\bold R\times [0,\infty)$. 
 Their union forms 
the ``backward light cone'' from the apex $(x,t)$ to the base  
line $\bold R\times\{0\}$.  
 For almost every such 
$y$, it must be that 
$$(\xi(s),s)\in \leb(f\circ\lambda)\cap \leb(\lambda)
$$
for almost every $s\in[0,t]$, because only a null set 
of points of $\bold R\times[0,\infty)$ fail to be 
Lebesgue points. Fix now such a $y$. 

The function $\gamma(s)=v(\xi(s),s)$, $0\le s\le t$, is 
Lipschitz, and hence $\gamma'(s)$ exists for a.e.\ $s$
and 
$$\gamma(t)-\gamma(0)=\int_0^t \gamma'(s)ds.
\tag 3.15
$$ 
Fix $s_0\in[0,t]$ such that all this holds: $s_0\in\bold T$, 
$\gamma'(s_0)$ exists, $(x_0,s_0), (\xi(s_0),s_0)
\in \leb(f\circ\lambda)\cap \leb(\lambda)$, and 
$s_0\in\leb\{f(\lambda(x_0,\cdot))\}
$. These conditions are satisfied by a.e.\ $s_0$. For such an $s_0$ we
prove that
$$\gamma'(s_0)=\frac{x-y}t \lambda(\xi(s_0),s_0)
-f( \lambda(\xi(s_0),s_0)).
\tag 3.16
$$
To prove (3.16) we calculate $\gamma'(s_0)$ as the limit 
$$\lim\Sb s\searrow s_0\\s\in\bold T\endSb 
\frac{\gamma(s)-\gamma(s_0)}{s-s_0}.
$$
So consider $s\in\bold T$, $s>s_0$. 
$$\aligned
&\frac{\gamma(s)-\gamma(s_0)}{s-s_0}
=\frac{v(\xi(s),s)-v(\xi(s_0),s_0)}{s-s_0}\\
=\;&\frac1{s-s_0}
\int_{\xi(s_0)}^{\xi(s)}\lambda(y,s)dy 
+\frac1{s-s_0}\left\{ \int_{x_0}^{\xi(s_0)}\lambda(y,s)dy 
- \int_{x_0}^{\xi(s_0)}\lambda(y,s_0)dy\right\}\\
&\ -\frac1{s-s_0}\int_{s_0}^s f(\lambda(x_0,\tau))d\tau\\
\equiv\;&A_1+A_2-A_3.
\endaligned
$$
Term by term: For $A_1$ first use $s-s_0=(\xi(s)-\xi(s_0))t/(x-y)$.
 By the assumption $s\in\bold T$ and Lemma 3.4(i), we can 
replace $s$ by an average over $[s_0,s]\cap\bold T$ at the 
expense of an $O(s-s_0)$ error term. Since $\bold T$ covers
almost all of $[s_0,s]$ we can ignore the restriction to $\bold T$
in the average. 
These steps give  
$$\aligned
A_1&=\frac1{s-s_0}
\int_{\xi(s_0)}^{\xi(s)}\lambda(y,s)dy \\
&=\frac{x-y}t \cdot \frac1{(s-s_0)(\xi(s)-\xi(s_0))}
\int_{s_0}^s d\tau \int_{\xi(s_0)}^{\xi(s)}\lambda(y,\tau)dy
+O(s-s_0) \\
&\longrightarrow \frac{x-y}t \lambda(\xi(s_0),s_0)
\endaligned
$$
as $s\searrow s_0$ along points of $\bold T$, by the 
assumption $(\xi(s_0),s_0)
\in \leb(\lambda)$. Next, 
$$A_2\longrightarrow  f( \lambda(x_0,s_0)) -f( \lambda(\xi(s_0),s_0))
$$
as $s\searrow s_0$ along points of $\bold T$,
by the assumptions $s_0\in\bold T$ and $(x_0,s_0), (\xi(s_0),s_0)
\in \leb(f\circ\lambda)$ and by Lemma 3.4(ii). Finally, 
$$A_3\longrightarrow  f( \lambda(x_0,s_0)) 
$$
by the assumption $s_0\in\leb\{f(\lambda(x_0,\cdot))\}
$. Altogether we have proved (3.16) for a.e.\ $s_0$. 

By (2.8),  (3.16) implies $\gamma'(s)\ge -g((x-y)/t)$ for a.e.\ $s$. Since
$v(y,0)=u_0(y)$, 
we can rewrite (3.15) as 
$$ v(x,t)=u_0(y)+\int_0^t \gamma'(s)ds\ge u_0(y)-tg((x-y)/t).
$$
This argument is valid for a.e.\ $y$, and consequently by the Hopf-Lax
formula (2.10), $v(x,t)\ge u(x,t)$. We have proved (3.14), and thereby
(3.6) for $x=x_0$. The choice of $x_0$ was such that it covers 
a.e.\ point in $\bold R$. 

For the converse, suppose 
$$\int_0^t f(\lambda(x,s))ds= \int_0^t f(\rho(x,s))ds
\qquad\text{for a.e.\ $(x,t)$.}
\tag 3.17
$$
 (3.8) is valid for both $\lambda$ and $\rho$
for a.e.\ $t$. Consider such a $t$ for which 
also (3.17) holds for a.e.\ $x$.
 The right-hand side of (3.8) is
$\int_{\bold R} dx\, \psi(x)\int_0^t f(\lambda(x,s))ds$, 
so by (3.17) 
its value is not changed by replacing $\lambda$ with $\rho$. 
We conclude that for a.e.\ $t$, 
$\int_{\bold R} \lambda(x,t) \psi(x)dx=
\int_{\bold R} \rho(x,t) \psi(x)dx$ for all 
test functions $\psi$. This implies $\lambda=\rho$ a.e.

%\enddocument

\hbox{}

\head 4. Second class particles in the variational coupling \endhead 

\subhead 4.1. Construction and variational coupling \endsubhead
First we  construct the process. 
Let $\{D_i:i\in\bold Z\}$ be a collection of mutually
independent rate 1 Poisson jump time processes on the time line 
$(0,\infty)$, defined on some probability
space $(\Omega,\Cal F, P)$. 
We  construct the process 
$\eta(t)=(\eta(i,t):i\in{\bold Z})$ on $\Omega$ in terms of 
``current particles.'' These form 
a process $z(t)= (z(i,t):i\in{\bold Z})$ 
 of labeled particles that move on  ${\bold Z}$. 
The location of the $i$th particle at time $t$ is 
$z(i,t)$, and these satisfy 
$$0\le z(i+1,t)-z(i,t)\le K\,.
\tag 4.1
$$
  In the graphical construction, 
 $z(i)$ attempts to jump 1 step to the {\it left} at 
epochs of $D_i$. If the jump 
violates (4.1), it is suppressed. 
 We can summarize the jump rule 
like this:
$$\aligned
&\text{If $t$ is an epoch of $D_i$, then}\\
&z(i,t)=\max\{ z(i,t-)-1, z(i-1,t-), z(i+1,t-)-K\}.
\endaligned
\tag 4.2
$$

To construct $z(\cdot)$ from rule (4.2)  
we exclude an expectional null set of realizations of 
$\{D_i\}$. For the remainder of this section,
fix a realization
$\{D_i\}$ with these properties:
\roster
\item "(4.3a)" There are no simultaneous jump attempts. 
\item "(4.3b)" Each
$D_i$ has only finitely may epochs in every bounded time
interval. 
\item "(4.3c)" Given any $t_1>0$, there are arbitrarily faraway
indices $i_0<<0<<i_1$ such that $D_{i_0}$ 
and $D_{i_1}$  have no epochs 
in the time interval $[0,t_1]$. 
\endroster
These properties are satisfied 
almost surely, so the evolution $z(t)$, $0\le t<\infty$,
 is well-defined
 almost everywhere on $\Omega$. Then the 
process $\eta(t)$ is defined by
$$\eta(i,t)=z(i,t)-z(i-1,t)\,.
\tag 4.4
$$ 
It should be clear that
$\eta(\cdot)$ satisfies the rules described 
in Section 2. 

Define a family $\{w^k:k\in{\bold Z}\}$ of auxiliary processes 
on $\Omega$. Each  
$w^k(t)=(w^k(i,t):i\in{\bold Z})$ is a process 
of the same type as $z(t)$, so  (4.1) is in force for all $w^k(i,t)$. 
The initial configuration of $w^k$ depends on the initial
position $z(k,0)$:
$$w^k(i,0)=\cases
z(k,0)\,, &i\ge 0\\
z(k,0)+Ki\,, &i<0. 
\endcases
%\tag 3
$$
The process $w^k$ reads its jump commands from the same 
Poisson processes $\{D_i\}$, but with a translated index:
at epochs $t$ of $D_{i+k}$,
$$w^k(i,t)=\max\{ w^k(i,t-)-1, w^k(i-1,t-), w^k(i+1,t-)-K\}.
$$

We prove a number of
pathwise properties for the processes
$z$ and $\{w^k\}$ that are valid under the fixed 
realization  $\{D_i\}$ that satisfies (4.3a--c). 
Since  (4.3a--c) are
valid almost surely, the properties we derive are 
almost sure properties, but  we leave the modifier ``almost
surely'' out
of the statements that follow. 

The coupling of the processes $z$ and $\{w^k\}$
 preserves ordering. One 
consequence is the following inequality, valid for
all $k<l$, all $i$ and all $t$:
$$w^l(i-l,t)\le w^k(i-k,t) + [z(l,0)-z(k,0)].
\tag 4.5
$$
The usefulness of the family
of processes $\{w^k\}$ lies in this ``variational coupling''
lemma:

\proclaim{Lemma 4.1} For all $i\in{\bold Z}$ and $t\ge 0$,
$$z(i,t)=\sup_{k\in{\bold Z}} w^k(i-k,t). 
\tag 4.6
$$
\endproclaim

In the case $K=\infty$ the supremum in (4.6) is 
over $\{k: k\le i\}$  because $w^k(i-k,t)=-\infty$ 
 for  $k>i$. 

For some purposes it is convenient to modify the processes 
$w^k$ so  that they
 start at level zero and advance in the
 positive direction. So define processes $\{\xi^k\}$ by
$$\xi^k(i,t)=z(k,0)-w^k(i,t).$$
The process 
$\xi^k$ does not depend on $z(k,0)$, and depends on 
the superscript $k$ only through a translation of the
$i$-index of the Poisson processes
$\{D_i\}$. 
 Initially 
$$\xi^k(i,0)=\cases
0\,, &i\ge 0\\
-Ki\,, &i<0\,.
\endcases
%\tag 4.9
$$
Dynamically, at epochs $t$ of $D_{i+k}$, 
$$
\xi^k(i,t)=\min\{ \xi^k(i,t-)+1, \xi^k(i-1,t-), \xi^k(i+1,t-)+K\}.
$$
These jumps preserve  the inequalities 
$$\xi^k(i,t)\le \xi^k(i-1,t)\qquad\text{and}\qquad
\xi^k(i,t)\le \xi^k(i+1,t)+K. 
%\tag 4.11
$$

Now (4.6) can be written as 
$$z(i,t)=\sup_{k\in{\bold Z}} \{z(k,0)-\xi^k(i-k,t)\}.
\tag 4.7
$$
This variational formula is a microscopic analogue of 
the Hopf-Lax formula (2.10), and  the basis for the proof
of the  hydrodynamic limit.  To prove Theorem 1, 
one shows that  the right-hand side of (4.7) converges to the 
right-hand side of (2.10). 
For more details about the construction, the variational
coupling lemma, and the proofs, the reader may consult 
 Sepp\"al\"ainen (1999).

\subhead 4.2. The second class particle \endsubhead
We investigate the position 
 $X(t)$  of a second class particle in the variational
coupling. Fix again a sample point so that we can regard
$\{D_i\}$, $(\eta(i,0))$, and $X(0)$ as deterministic, 
with (4.3a--c) satisfied by $\{D_i\}$. 

As defined in Section 2.2, $X(t)$ is the location 
of the unique discrepancy between two processes 
$\eta$ and $\etatil$ that satisfy initially 
$\etatil(X(0),0)=\eta(X(0),0)+1$ and $\etatil(i,0)=\eta(i,0)$
for $i\ne X(0)$. 
In particular, the second class particle
is not counted in the $\eta$-variables and by the 
$K$-exclusion rule  always
 $\eta(X(t),t)\le K-1$. 
 
We define $\eta$ and $\etatil$ by (4.4),  in terms
of processes $z$ and $\ztil$ that initially satisfy
$$
\text{$\ztil(i,0)=z(i,0)$ for $i\le X(0)-1$
and $\ztil(i,0)=z(i,0)+1$ for $i\ge X(0)$. }
$$
Processes $z$ and $\ztil$ follow the same Poisson
processes $\{D_i\}$, so this is the basic coupling. 
The location $X(t)$ is always uniquely defined 
as the location of the discrepancy, as shown in 
this lemma:

\proclaim{Lemma 4.2} For all $t\ge 0$ there is a finite
location $X(t)$ such that 
$$
\text{$\ztil(i,t)=z(i,t)$ for $i\le X(t)-1$
and $\ztil(i,t)=z(i,t)+1$ for $i\ge X(t)$. }
\tag  4.8
$$
\endproclaim

\demo{Proof} Given any finite time horizon $t_0$,
find indices $\cdots<i_{-1}<i_0<i_1<\cdots$
such that $D_{i_m}$ has no jump epochs in the 
time interval $(0,t_0]$ for all $m$,
 $i_m\to\pm\infty$ as $m\to\pm\infty$,
and $i_0<X(0)<i_1$.
Particles $z(i_m)$ and $\ztil(i_m)$ do not 
move up to time $t_0$. In each finite segment 
$(i_{m-1},i_m)$ there are only finitely many
jump epochs up to time $t_0$. Arrange these finitely many
jump epochs in temporal order, and 
then  prove by induction that 
the lemma  holds after each jump time. 
\qed
\enddemo

(4.8) defines $X(t)$ in our framework. 
Next we give $X(t)$ a representation that is a microscopic
analogue of the characteristics.

\proclaim{Proposition 4.1} 
We have the formulas 
$$X(t)=\inf\{ i\in{\bold Z}: \text{$z(i,t)=z(k,0)-\xi^k(i-k,t) $ for
some $k\ge X(0)$}\} 
\tag 4.9
$$
and 
$$X(t)=1+\sup\{ i\in{\bold Z}: \text{$\ztil(i,t)=\ztil(k,0)-
\xi^k(i-k,t) $ for
some $k< X(0)$}\}.
\tag 4.10
$$
\endproclaim 

\demo{Proof} Consider these four statements:
$$
\text{If $i<X(t)$, then $z(i,t)>z(k,0)-\xi^k(i-k,t) $ for
all $k\ge X(0)$.}
\tag  4.11
$$
$$
\text{If $i\ge X(t)$, then $z(i,t)=z(k,0)-\xi^k(i-k,t) $ for
some $k\ge X(0)$.}
\tag 4.12
$$
$$
\text{If $i<X(t)$, then $\ztil(i,t)=\ztil(k,0)-\xi^k(i-k,t) $ for
some $k< X(0)$.}
\tag 4.13
$$
$$
\text{If $i\ge X(t)$, then $\ztil(i,t)>\ztil(k,0)-\xi^k(i-k,t) $ for
all $k< X(0)$.}
\tag 4.14
$$
These statements imply that the infimum and 
supremum in the formulas (4.9--10) are well-defined and
finite, and also prove the formulas themselves. 

To contradict (4.11), suppose  $i<X(t)$ and
$$z(i,t)=z(k,0)-\xi^k(i-k,t)
\tag 4.15
$$ 
for
some $k\ge X(0)$. Then by (4.8)
$$\ztil(i,t)=\ztil(k,0)-\xi^k(i-k,t)-1$$
which contradicts the variational formula (4.7) 
for process $\ztil$. This contradiction proves (4.11), and also
implies that there must exist some $k< X(0)$
such that (4.15) holds. By (4.8), (4.15) then turns into  
$$\ztil(i,t)=\ztil(k,0)-\xi^k(i-k,t), $$
and this proves (4.13). 

(4.12) and (4.14) are proved by a similar argument that starts
by contradicting (4.14). 
\qed
\enddemo

The second class particle may jump either left or right, 
but a certain monotonicity can be found from the 
representation (4.9): 

\proclaim{Proposition 4.2} 
Let 
$$k_1(t)=\sup\{k\ge X(0): z(X(t),t)= w^k(X(t)-k,t)\}
$$
denote the maximal $k$ 
that satisfies the requirement in
formula {\rm (4.9)}, or infinity. Then $k_1(t)$ is
 nondecreasing 
as a function of time. 
  If there are initially arbitrarily 
large indices $j$ such that $\eta(j,0)\le K-1$, 
then $k_1(t)$ is  finite for all $t$. 
\endproclaim 

We shall not prove Proposition 4.2 for we make no use of 
it in this paper.

\hbox{}

\head 5. Proofs of the laws of large numbers \endhead

Now assume we are in the setting described in the 
paragraph preceding Theorem 1. The initial particle 
configurations $(\eta_n(i,0):i\in\bold Z)$, $n\in\bold N$,
 and 
the Poisson jump time processes $\{D_i\}$ are 
defined on a probability space $(\Omega, \Cal F, P)$.
Define initial 
configurations $(z_n(i,0):i\in\bold Z)$ by setting 
 $z_n(0,0)=0$ and using (4.4). Construct the processes
$z_n(\cdot)$ as in Section 4.1, and then 
define the processes $\eta_n(\cdot)$
by (4.4) in terms of
$z_n(\cdot)$. 
The location of the second class particle 
is given by  the variational representation 
in Proposition 4.1: 
for the $n$th process
$$
\aligned
X_n(t)&=\inf\{ i\in{\bold Z}: \text{$z_n(i,t)=w_n^k(i-k,t) $ for
some $k\ge X_n(0)$}\}\\
 &=\inf\{ i\in{\bold Z}: \text{$z_n(i,t)=z_n(k,0)-\xi^k(i-k,t) $}\\
 &\quad\qquad\qquad\qquad \text{ for
some $k\ge X_n(0)$}\}.
\endaligned
%\tag 
$$

Theorems 3 and 4 rely on the basic 
hydrodynamic limits, which we summarize here
from \Se\  in terms of the 
$z$-variables. The following statements all hold
almost surely:
$$\lim_{n\to\infty}\frac1n z_n([ny],0)=u_0(y)
\quad\text{for all $y\in\bold R$;}
\tag 5.1
$$
$$\lim_{n\to\infty}\frac1n \xi^{[nr]}([nx],nt)=tg(x/t) 
\quad\text{for all $x,r\in\bold R$ and $t>0$;}
\tag 5.2
$$
and 
$$\lim_{n\to\infty}\frac1n z_n([nx],nt)=u(x,t)
\quad\text{for all $x\in\bold R$ and $t>0$.}
\tag 5.3
$$
The function $u_0$ is defined by (2.7) from the given initial profile 
$\rho_0$. 
 (5.1) is a restatement of assumption  (2.3). 
(5.2) is proved by subadditivity, and 
this limit defines $g$. 
Then $u(x,t)$ is defined in terms of $u_0$
and $g$ by the Hopf-Lax formula (2.10). Next
(5.3) is proved, by (5.1), (5.2), and the variational
coupling (4.7).  Lastly, 
 $f$ is defined as the conjugate of $g$,
and $u(x,t)$ is identified as the viscosity solution
of (2.11).

The following estimate  [Sepp\"al\"ainen (1999), Lemma 6.1]
reduces the range of indices
needed in the variational formula. It will be used 
several times in the proofs. 
\proclaim{Lemma 5.1} 
Define 
$$\zeta_n(i,l,t)=\max_{i-l\le k\le i+l}
\{ z_n(k,0)-\xi^k(i-k,t)\}.
$$
Fix $r>t$. Then 
 there exists
a finite positive constant  $C=C(r,t)$ such that
for all $i\in\bold Z$ and all $n\ge 1$, 
$$
P\bigl( z_n(i,nt)\ne \zeta_n(i,nr,nt) \bigr)\le e^{-Cn}.
\tag 5.4
$$
\endproclaim

First we prove Theorem 3, so recall assumption (2.20):
$$\lim_{n\to\infty}\frac1n X_n(0)=b \quad\text{a.s.}
\tag 5.5
$$

\proclaim{Proposition 5.1} Assume {\rm (2.3)} and {\rm (2.20)}. 
Then 
$$\limsup_{n\to\infty} n^{-1}X_n(nt)\le x^+(b,t)
\quad\text{ almost surely.}
$$
\endproclaim

\demo{Proof} 
 Fix $x_1>x^+(b,t)$. 
We shall show that 
$\limsup_{n\to\infty} n^{-1}X_n(nt)\le x_1$ almost surely.
By definition $x^+(b,t)=\sup\{x: y^-(x,t)\le b\},$
so 
$y^-(x_1,t)>b$. Choose $c_1$ so that 
$b<c_1<y^-(x_1,t)$. 

Let 
$\Gamma$ be the event 
 on which 
$$\text{$X_n(nt)>nx_1$ for infinitely many $n$,}
\tag 5.6
$$
assumptions (4.3a--c) of Section 4.1 hold so that we can use
variational formulas (4.7) and (4.9), 
the laws of large numbers (5.1)--(5.3) and (5.5) hold, and 
$$z_n([nx_1],nt)=w_n^k([nx_1]-k,nt)
\tag 5.7
$$ 
for some  $k\ge [nc_0]$ 
for all large enough $n$, 
 for some fixed finite $c_0$. Except for (5.6) the conditions
of $\Gamma$ hold with probability one, so it suffices to show 
that $\Gamma$ is empty. 
 
Suppose
$\Gamma$ is not empty. 
Fix a sample point in $\Gamma$, and suppose $n$ is such that
$X_n(nt)>nx_1$. Then 
$z_n([nx_1],nt)$ $\ne$ $ w_n^k([nx_1]-k,nt)$ for all 
$k\ge X_n(0)$, so there must exist a $k<X_n(0)$ such that
(5.7) holds. 
Since limit (5.5) holds on $\Gamma$,
  $X_n(0)\le [nc_1]$ for large enough $n$. Thus 
 for a subsequence of $n$'s there exists a 
$k$ in (5.7) such that $[nc_0]\le k\le [nc_1]$. 

Pick a partition 
$c_0=y_0<y_1<\cdots<y_m=c_1$ of $[c_0,c_1]$ with mesh
$\Delta y=\max(y_{i+1}-y_i)$. Fix $i$ so that 
for infinitely many $n$ in the subsequence picked above, 
some $k$ in (5.7) satisfies 
$[ny_i]\le k \le [ny_{i+1}]$. By (4.5) we have, for 
infinitely many $n$, 
$$
\aligned
z_n([nx_1],nt)&=w_n^k([nx_1]-k,nt)\\
&\le w_n^{[ny_i]}([nx_1]-[ny_i],nt)+\bigl\{ z_n([ny_{i+1}],0)-
z_n([ny_{i}],0)\bigr\}.
\endaligned
\tag 5.8
$$
Divide by $n$ and let $n\to\infty$ along the 
appropriate subsequence of $n$'s. 
 Since the laws of
large numbers hold for our fixed sample point, 
we get 
$$u(x_1,t)\le u_0(y_i)-tg\biggl(\frac{x_1-y_i}{t}\biggr)
+\bigl\{ u_0(y_{i+1})-
u_0(y_{i})\bigr\}.
\tag 5.9
$$
Take the mesh $\Delta y\to 0$. By continuity and
compactness, we conclude that the Hopf-Lax formula
(2.10) for $u(x_1,t)$ has a maximizer in $[c_0,c_1]$, 
 so that 
$y^-(x_1,t)\le c_1$.
 This contradiction with the choice of $c_1$
 implies that $\Gamma$ is empty.
\qed
\enddemo

\proclaim{Proposition 5.2} Assume {\rm (2.3)} and {\rm (2.20)}. Then 
$$\liminf_{n\to\infty} n^{-1}X_n(nt)\ge x^-(b,t)
\qquad\text{ almost surely.}
$$
\endproclaim

\demo{Proof} Fix $x_1<x^-(b,t)$. Consider the event $\Gamma$
on which  $X_n(nt)<nx_1$ for infinitely many $n$,
assumptions (4.3a--c) hold, and 
 the laws of large numbers (5.1)--(5.3), and (5.5) hold.
 We assume that 
$\Gamma$ has positive probability, and  derive a contradiction. 

(4.12) and $X_n(nt)<nx_1$ imply that for some
 $k\ge X_n(0)$, 
$$z_n([nx_1],nt)= z_n(k,0)-\xi^k([nx_1]-k,nt).
\tag 5.10
$$
As in Section 4.2, we have the processes $\ztil_n$,
coupled with $z_n$ through the Poisson processes $\{D_i\}$,
that satisfy
$$
\aligned
&\text{$\ztil_n(i,nt)=z_n(i,nt)$ for $i<X_n(nt)$}\\
\text{and }\;
&\text{$\ztil_n(i,nt)=z_n(i,nt)+1$ for $i\ge X_n(nt)$.}
\endaligned
\tag 5.11
$$
(4.10) implies that 
$$\ztil_n([nx_1],nt)= \ztil_n(k,0)-\xi^k([nx_1]-k,nt).
\tag 5.12
$$
 By (4.14), (5.12) cannot hold for any
$k<X_n(0)$. 
By applying Lemma 5.1 to the processes $\ztil_n$ and by
shrinking $\Gamma$ by no more than a null set, 
we may assume that for large enough $n$ the $k$ in (5.12)
must lie in the range $X_n(0)\le k\le [nc_1]$ for a
constant $c_1$. Then by (5.5) we conclude  
that for large $n$, $[nc_0]\le k\le [nc_1]$ 
where  $c_0$ can be taken arbitrarily close to $b$. 
By left-continuity of $x^-(\cdot,t)$ and the assumption
$x_1<x^-(b,t)$, we may  choose $c_0$ so that $x_1<x^-(c_0,t)$.

Now starting from (5.12) repeat 
the argument around (5.8)--(5.9) for $\ztil_n$. 
By (5.11) $\ztil_n$ satisfies the limits (5.1) and  (5.3). 
As in the step after (5.9), we get that $u(x_1,t)$ has 
a maximizer in $[c_0,c_1]$ so that $y^+(x_1,t)\ge c_0$. 
By definition (2.18) this implies $x^-(c_0,t)\le x_1$
and contradicts the choice of $c_0$. 
\qed
\enddemo

Propositions 5.1--5.2  prove Theorem 3.

\demo{Proof of Theorem 4(a)} 
Assume 
$x(y,t)=x^\pm(y,t)$ for all $y\in\bold R$. 
Let $\Gamma$ be the event on which (4.3a--c) and (5.1)--(5.3)
hold, and 
also for
 large enough $n$ 
$$-nb_0\le X_n(0)\le nb_0 
\tag 5.13
$$
for some constant $b_0$.
We prove that 
$$
\lim_{n\to\infty} \left| \frac{X_n(nt)}n
- x\left( \frac{X_n(0)}n \,,\,t\right) \right|=0.
\tag 5.14
$$
 almost
surely on $\Gamma$. By assumption (2.23), Theorem 4(a)
follows by repeating this for a sequence of $b_0$'s 
increasing to $\infty$. 

Since $x^-(y, t)$ is left-continuous and 
$x^+(y, t)$  right-continuous in $y$, the assumed
unique characteristic $x(y,t)=x^\pm(y,t)$ 
is continuous in $y$. Given $\e>0$, choose
a partition $-b_0=r_0<r_1<\cdots<r_h=b_0$
so that 
$x(r_{j+1},t)<x(r_j,t)+\e/2$. Set $x_j=x(r_j,t)+\e$.
Shrink $\Gamma$ by a null set to guarantee that 
a.s.\ for every $j=0,\ldots,h$, for large enough
$n$ there exists a $k\ge nc_0$
such that 
$$z_n([nx_j],nt)=w_n^k([nx_j]-k,nt).
$$
This can be done by Lemma 5.1 if $c_0$ is chosen small enough. 

To contradict (5.14), suppose $\Gamma$ contains a sample point
for which 
$$X_n(nt)> nx(n^{-1}X_n(0),t)+n\e
\tag 5.15
$$
for infinitely many $n$. Restrict this subsequence
of $n$'s further so that for some fixed $j$, 
$$nr_j\le X_n(0)\le nr_{j+1}
\tag 5.16
$$ 
for all these $n$. By the monotonicity of $x(\cdot, t)$,
(5.15--16) imply that  
$X_n(nt)> nx_j$ 
for infinitely many $n$. For each such $n$  there exists
 $k<X_n(0)$ and no $k\ge X_n(0)$ such that
$z_n([nx_j],nt)$ $=$ $w_n^k([nx_j]-k,nt).
$
By the assumptions made above, 
$nc_0\le k\le nr_{j+1}$. 
 Repeat the partitioning
and limiting argument around  (5.8)--(5.9) 
to conclude  that $y^-(x_j,t)\le r_{j+1}$, and then by  
definition  $x^+(r_{j+1},t)\ge x_j$. Now 
we have contradicted 
$x_j=x(r_j,t)+\e>x(r_{j+1},t)+\e/2$.

It remains to show that assuming 
$X_n(nt)< nx(n^{-1}X_n(0),t)-n\e$
for infinitely many $n$ also leads to a contradiction.
Combine this with (5.16) to get 
$$X_n(nt)< nx(n^{-1}X_n(0),t)-n\e\le 
nx(r_{j+1},t)-n\e\equiv n x'_j.
$$
As in the proof of Proposition 5.2, switch to the 
$\ztil_n$ processes to get  
$$\ztil_n([nx'_j],nt)=\ztil_n(k,0)-\xi^k([nx'_j]-k,nt)$$
for some $k$ such that  $nr_j\le X_n(0)\le k\le nc_1$,
where $c_1$ is a fixed constant. 
From this  conclude that 
$y^+(x'_j,t)\ge r_j$, which implies 
$x(r_j,t)=x^-(r_j,t)\le x'_j=x(r_{j+1},t)-\e$, 
a contradiction. 
\qed
\enddemo

To prove Theorem 4(b), first a lemma.

\proclaim{Lemma 5.2} Let $c$ be a finite constant.

 {\rm (a)} Suppose that 
 $X_n(0)\ge nc$ for large enough $n$, almost surely. 
Then for any $x<c-t$, 
$X_n(nt)\ge nx$ for large enough $n$, almost surely. 

 {\rm (b)} Suppose that 
 $X_n(0)\le nc$ for large enough $n$, almost surely. 
Then for any $x>c+t$, 
$X_n(nt)\le nx$ for large enough $n$, almost surely. 
\endproclaim

\demo{Proof} (a)  
 Pick $\e>0$ so that $x<c-t-2\e$,
and set $c_1=c-\e$. 

As an intermediate claim, we prove that
almost surely,  for large enough
$n$, 
$$\xi^k([nx]-k,nt)=-K\cdot([nx]-k)\quad
\text{for all $k\ge [nx]+[n(t+\e)]$.} 
\tag 5.17
$$
In other words, the claim is that $\xi^k([nx]-k,\cdot)$ has
not jumped by time $nt$. The 
 time of the first jump of $\xi^k([nx]-k,\cdot)$
is distributed as $\sum_{j=0}^{k-[nx]}Y_j$ where
$\{Y_j\}$ are i.i.d.\ rate 1 exponential random variables. 
 The estimate 
$P(\sum_{j=1}^N Y_j\le Ns)\le \exp[-N\kappa(s)]$ is valid for 
 $0<s<1$. The  rate function is 
$\kappa(s)=s-1-\log s$ for $s>0$, which satisfies 
$\kappa'(s)<0<\kappa(s)$ 
for $0<s<1$.
Thus 
$$\aligned
&P\left\{\text{(5.17) fails}\right\}
\le \sum_{k\ge [nx]+[n(t+\e)]} 
P\left( \sum_{j=0}^{k-[nx]}Y_j \le nt\right)\\
&\le \sum_{i\ge [n(t+\e)]} 
\exp\left[-i\kappa\left({nt}/{i}\right)\right]
\le C_1e^{-C_2n}
\endaligned
$$
for finite positive constants $C_1,C_2$. This proves
the claim, and (5.17) holds for large enough $n$ a.s.

Since $\eta_n(X(0),0)\le K-1$ and $[nc_1]<X(0)$, 
$$z_n(k,0)-z_n([nc_1],0)=\sum_{i=[nc_1]+1}^k\eta_n(i,0)
\le K(k-[nc_1])-1$$
 for any $k\ge X_n(0)$. For the fixed $x$ any 
 $k\ge X_n(0)$ falls within the range in (5.17) for
large enough $n$.  So for large
enough $n$, 
$$\aligned
&z_n(k,0)-\xi^k([nx]-k,nt)=z_n(k,0)+K([nx]-k)\\
&\le z_n([nc_1],0)+K([nx]-[nc_1])-1 
  =z_n([nc_1],0)-\xi^{[nc_1]}([nx]-[nc_1],nt)-1.
\endaligned
$$
The last equality comes from (5.17) again because 
$[nc_1]$ also falls within the range of $k$'s in (5.17) for
large enough $n$. 
The string of inequalities shows that 
 no $k\ge X_n(0)$ can give the supremum  in the microscopic
variational formula (4.7) for $z_n([nx],nt)$, because each $k\ge X_n(0)$ is
inferior to $[nc_1]$. 
But if $X_n(nt)<[nx]$, by (4.12) some 
$k\ge X_n(0)$ would be  a maximizer. 
This contradiction shows that $X_n(nt)\ge[nx]$ for
large enough $n$, a.s.

The proof of part (b) is similar, but instead of (5.17)
use that $\xi^k([nx]-k,nt)=0$ for all $k\le [nx]-[n(t+\e)]$. 
\qed
\enddemo

\demo{Proof of  Theorem 4(b)}
The goal is to show that 
$$
\lim_{n\to\infty} \left| \frac{X_n(nt)}n
- x\left( \frac{X_n(ns)}n \,;s,\,t\right) \right|=0
\tag 5.18
$$

Let $\ubar_0(x)=u(x,s)$,
$\ubar(x,\tau)=u(x,s+\tau)$ for $\tau>0$, and let
$\xbar(b;0,\tau)$ denote characteristics for the 
solution $\ubar$. Then 
 $x^\pm(b;s,t)=\xbar^\pm(b;0,t-s)$.
From the assumptions and Proposition 3.1 
we know $\ubar$ has unique characteristics:
$\xbar(b;0,\tau)=\xbar^\pm(b;0,\tau)$. 

Let $\zbar_n(i,t)=z_n(i,ns+t)$ be the $n$th process 
restarted at time $ns$. 
Let $\Xbar_n(t)$ be the location
of a second class particle in the process $\zbar_n$, 
with initial location $\Xbar_n(0)=X_n(ns)$. 
In restarted terms, (5.18) is the same as 
$$
\lim_{n\to\infty} \left| \frac{\Xbar_n(n(t-s))}n
- \xbar\left( \frac{\Xbar_n(0)}n \,;0,\,t-s\right) \right|=0
\ \ \ \text{a.s.}
\tag 5.19
$$
 (5.19) follows from Theorem 4(a) if we 
 check that
the ingredients of  the proof are valid for 
the processes $\zbar_n$ and $\Xbar_n$. 

\hbox{}

(i) We can obtain the variational formulas (4.7) and 
(4.9) for $\zbar_n$ and $\Xbar_n$ in a form where 
the $\xi$-part depends on $n$ through the time shift. 
For any fixed restarting time $\tau\ge 0$, 
let $\theta_\tau$ denote the
restarting  
operation on the probability space
of the Poisson processes $\{D_i\}$:
if the epochs of $D_i$ in $(0,\infty)$ are
$t_1<t_2<t_3<\cdots$ with $t_{j-1}\le \tau <t_j$, then 
the epochs of $\theta_\tau D_i$ in $(0,\infty)$ are 
$t_j-\tau<t_{j+1}-\tau<t_{j+2}-\tau<\cdots$. 

The same proofs that originally gave Lemma 4.1 and 
Proposition 4.1 also yield, 
 for $t>\tau$, 
$$z(i,t)=\sup_k\{z(k,\tau)-\xi^{k}(i-k,t-\tau)\circ\theta_\tau\} 
\tag 5.20
$$
and 
$$X(t)=\inf\{i: 
\text{$z(i,t)=z(k,\tau)-\xi^{k}(i-k,t-\tau)\circ\theta_\tau$ for some
$k\ge X(\tau)$} \}.
\tag 5.21
$$

Now for process $z_n$ take $\tau=ns$, and (5.20)--(5.21) give 
for $t>0$
$$\zbar_n(i,t)
=\sup_k\{\zbar_n(k,0)-\xi^{k}(i-k,t)\circ\theta_{ns}\} 
%\tag 
$$
and 
$$\Xbar_n(t)=\inf\{i: 
\text{$\zbar_n(i,t)=\zbar_n(k,0)-
\xi^{k}(i-k,t)\circ\theta_{ns}$ for some
$k\ge X_n(ns)$} \}.
%\tag 
$$

\hbox{}

(ii) Limits (5.1) and (5.3) are true by definition of $\zbar_n$ and 
$\ubar$. Limits (5.2) for $\xi^{[nr]}([nx],nt)\circ\theta_{ns}$ are true
by the proofs in \Se\ because these limits are
derived from summable deviation estimates and the 
Borel-Cantelli lemma. The distribution 
of $\xi^{[nr]}$ is not changed by the time shift,
 so the same deviation estimates are valid. 

\hbox{}

(iii) Finally assumption (2.25) for $\Xbar_n(0)=X_n(ns)$ is  a
consequence of (2.25) for $X_n(0)$ and Lemma 5.2. 
\qed
\enddemo

\hbox{}

\head References \endhead

\hbox{}

\flushpar
Dafermos, C. M. (1977). Generalized characteristics and 
the structure of solutions of hyperbolic 
conservation laws. Indiana Univ. Math. J. 26 1097--1119.

\hbox{}

\flushpar  
L. C. Evans (1998). Partial Differential Equations. 
American Mathematical Society.

\hbox{}

\flushpar  
Ferrari, P. A. (1992). Shock fluctuations in asymmetric
simple exclusion. \ptrf\ 91 81--101. 

\hbox{}

\flushpar  
P. A. Ferrari and L. R. G. Fontes (1994a). 
Shock fluctuations in the asymmetric
simple exclusion process. \ptrf\ 99 305--319. 

\hbox{}

\flushpar  
P. A. Ferrari and L. R. G. Fontes (1994b). Current
fluctuations for the asymmetric simple exclusion
process. \ap \ 22 820--832. 

\hbox{}

\flushpar  
P. A. Ferrari, L. R. G. Fontes and Y. Kohayakawa (1994).
Invariant measures for a two-species asymmetric process.
\jsp\ 76 1153--1177.

\hbox{}

\flushpar  
Ferrari, P. A. and Kipnis C. (1995). Second class particles
in the rarefaction fan. \aihp\ 31 143--154. 

\hbox{}

\flushpar  
H. Ishii (1984). Uniqueness of unbounded viscosity solutions
of Hamilton-Jacobi equations. Indiana Univ. Math. J. 33 721--748.
  
\hbox{}

\flushpar  
C. Kipnis and C. Landim (1999). Scaling Limits 
of Interacting Particle Systems.
 Grundlehren der mathematischen Wissenschaften, vol 320,
Springer Verlag, Berlin.

\hbox{}

\flushpar  
T. M. Liggett (1999). 
 Stochastic Interacting Systems: Contact, Voter
and Exclusion Processes. 
Springer-Verlag,  New York. 
 
\hbox{}

\flushpar  
Rezakhanlou, F. (1995). Microscopic structure of shocks
in one conservation laws. \aihpnl\ 12 119-153. 

\hbox{}

\flushpar  
R. T. Rockafellar (1970).  
Convex Analysis. Princeton University Press.          

\hbox{}

\flushpar  
 Rost, H. (1981). 
Non-equilibrium behaviour of a many particle
process: Density profile and local equilibrium. 
\zwvg\  58 41--53. 

\hbox{}

\flushpar  
T. Sepp\"al\"ainen (1996). A microscopic
model for the Burgers equation and longest increasing
subsequences.   Electronic J.
Probab.  1,  Paper 5, 1--51.

\hbox{}

\flushpar  
T. Sepp\"al\"ainen (1998a). Hydrodynamic scaling, convex
duality, and asymptotic shapes of growth models. 
 \mprf \ 4 1--26. 

\hbox{}

\flushpar  
T. Sepp\"al\"ainen (1998b). 
Coupling the totally asymmetric simple exclusion process
with a moving interface.
Proceedings of I Escola Brasileira de Probabilidade 
(IMPA, Rio de Janeiro, 1997). 
 \mprf \ 4 593--628.

\hbox{}

\flushpar  
T.\ Sepp\"al\"ainen (1999). Existence of hydrodynamics
for the totally asymmetric
simple $K$-exclusion process. 
 Ann.\ Probab. 27 361--415.

\hbox{}

\flushpar  
T.  Sepp\"al\"ainen (2000).
 A variational coupling for a totally asymmetric
 exclusion process with long jumps but no passing.
To appear in the Proceedings
of the Workshop on Hydrodynamic Limits, The Fields Institute,
Toronto, 1998. 

\enddocument

\hbox{}

\flushpar
Department of Mathematics\hfil\break
Iowa State University\hfil\break
Ames, Iowa 50011\hfil\break
E-mail: seppalai\@iastate.edu\hfil\break

\enddocument